\pdfoutput=1
\RequirePackage{ifpdf}
\ifpdf 
\documentclass[pdftex]{sigma}
\else
\documentclass{sigma}
\fi

\usepackage{enumerate}
\usepackage[all]{xy}

\numberwithin{equation}{section}

\newtheorem{Lemma}{Lemma}[section]
\newtheorem{Theorem}[Lemma]{Theorem}
\newtheorem{Proposition}[Lemma]{Proposition}
\newtheorem{Corollary}[Lemma]{Corollary}

{\theoremstyle{definition}
\newtheorem{Definition}[Lemma]{Definition}
\newtheorem{Example}[Lemma]{Example}
\newtheorem{Remark}[Lemma]{Remark}
\newtheorem{Notation}[Lemma]{Notation}
}

\begin{document}

\allowdisplaybreaks

\renewcommand{\thefootnote}{$\star$}

\newcommand{\arXivNumber}{0907.2891}

\renewcommand{\PaperNumber}{055}

\FirstPageHeading

\ShortArticleName{Non-Compact Symplectic Toric Manifolds}

\ArticleName{Non-Compact Symplectic Toric Manifolds\footnote{This paper is a~contribution to the Special Issue
on Poisson Geometry in Mathematics and Physics.
The full collection is available at \href{http://www.emis.de/journals/SIGMA/Poisson2014.html}
{http://www.emis.de/journals/SIGMA/Poisson2014.html}}}

\Author{Yael KARSHON~$^\dag$ and Eugene LERMAN~$^\ddag$}

\AuthorNameForHeading{Y.~Karshon and E.~Lerman}

\Address{$^\dag$~Department of Mathematics, University of Toronto,\\
\hphantom{$^\dag$}~40~St.~George Street, Toronto, Ontario, Canada M5S 2E4}
\EmailD{\href{mailto:karshon@math.toronto.edu}{karshon@math.toronto.edu}}

\Address{$^\ddag$~Department of Mathematics, The University of Illinois
at Urbana-Champaign,\\
\hphantom{$^\ddag$}~1409 W.~Green Street, Urbana, IL 61801, USA}
\EmailD{\href{mailto:lerman@math.uiuc.edu}{lerman@math.uiuc.edu}}

\ArticleDates{Received August 15, 2014, in f\/inal form July 10, 2015; Published online July 22, 2015}

\Abstract{A key \looseness=-1 result in equivariant symplectic geometry
is Delzant's classif\/ication of compact connected symplectic toric manifolds.
The moment map induces an embedding of
the quotient of the manifold by the torus action into the dual of the
Lie algebra of the torus; its image is a unimodular
(``Delzant'') polytope; this gives a bijection between unimodular polytopes
and isomorphism classes of compact connected symplectic toric manifolds.
In this paper we extend Delzant's classif\/ication to non-compact
symplectic toric manifolds.  For a~non-compact symplectic toric
manifold the image of the moment map need not be convex and the
induced map on the quotient need not be an embedding.  Moreover, even
when the map on the quotient is an embedding, its image no longer
determines the symplectic toric manifold; a degree two characteristic
class
on the quotient makes an appearance.  Nevertheless, the quotient is a
manifold with corners, and the induced map from the quotient to the
dual of the Lie algebra is what we call a unimodular local embedding.
We classify non-compact symplectic toric manifolds in terms of
manifolds with corners equipped with degree two cohomology classes
and
unimodular local embeddings into the dual of the Lie algebra
of the correspon\-ding torus. The main new ingredient is the construction
of a symplectic toric manifold from such data.
The proof passes through an equivalence of categories
between symplectic toric manifolds and symplectic toric bundles
over a f\/ixed unimodular local embedding.
This equivalence also gives a geometric interpretation
of the degree two cohomology class.}


\Classification{53D20;  53035; 14M25; 37J35}

\renewcommand{\thefootnote}{\arabic{footnote}}
\setcounter{footnote}{0}

\section{Introduction}

In the late 1980s, Delzant classif\/ied compact connected symplectic
toric manifolds \cite{De} by showing that the map
\begin{gather*}
\text{symplectic toric manifold} \ \mapsto \ \text{its moment map image}
\end{gather*}
is a bijection onto the set of unimodular (also referred to
as ``smooth'' or ``Delzant'') polytopes.
This beautiful work has been widely inf\/luential.
The goal of this paper is to extend Delzant's classif\/ication
theorem to non-compact manifolds.

Delzant's classif\/ication is built upon convexity and connectedness
theorems of Atiyah and Guillemin--Sternberg \cite{At,GS:convexity}.
Compactness plays a crucial role in the proof of these theorems.
Indeed, for a non-compact symplectic toric manifold the moment map
image need not be convex and the f\/ibers of the moment map need not be
connected. And even when the f\/ibers of the moment map are connected
the moment map image need not uniquely determine the corresponding
symplectic toric manifold.  Thus, the passage to noncompact symplectic
toric manifolds requires a dif\/ferent approach.  As a f\/irst step we
make the following observation (the proof is given in
Appendix~\ref{sec:background}):

\begin{Proposition} \label{prop:unimod-emb 1}
  Let $(M,\omega, \mu)$ be a symplectic toric $G$-manifold.
  Then the quotient $M/G$ is naturally a manifold with corners
  and the induced map
\begin{gather*} 
  \bar{\mu} \colon \ M/G\to {\mathfrak g}^*, \qquad \bar{\mu}(G\cdot x) := \mu(x)
\end{gather*}
is a unimodular local embedding.
$($See Definitions~{\rm \ref{quotient}} and~{\rm \ref{def:unimodular emb}.)}
\end{Proposition}

\begin{Definition} \label{def:orbital_moment_map}\sloppy
  Given a symplectic toric $G$-manifold $(M,\omega, \mu)$
  and a $G$-quotient map \mbox{$\pi \colon M \to W$},
  we refer to the map $\psi\colon W \to {\mathfrak g}^*$
  that is def\/ined by $\mu = \psi \circ \pi$ as the {\em orbital moment map}.
\end{Definition}

See Remarks~\ref{add_rk:1} and \ref{add_rk:2} for the origin of the notion of
orbital moment map and its relation to  developing map in
af\/f\/ine geometry.  The fact that the quotient~$M/G$ is a manifold with
corners is closely related to the fact that for a completely
integrable system with elliptic singularities the space of tori is a~manifold with corners~\cite{BM,Y}.

Our classif\/ication result can be stated as follows.
\begin{Theorem} \label{thm:main}
Let ${\mathfrak g}^*$ be the dual
  of the Lie algebra of a torus $G$ and $\psi\colon W\to {\mathfrak g}^*$ a unimodular
  local embedding of a manifold with corners.  Then
\begin{enumerate}[$1.$]\itemsep=0pt
\item\label{1}
  There exists a symplectic toric $G$-manifold $(M,\omega, \mu)$
  with $G$-quotient map $\pi\colon M \to W$
  and orbital moment map $\psi$.
\item \label{2} The set of isomorphism classes
  of symplectic toric $G$-manifolds $M$
  with $G$-quotient map $\pi\colon M \to W$
  and orbital moment map $\psi$ is in bijective
  correspondence with the set of cohomology classes
\begin{gather*}
H^2 (W, {\mathbb Z}_G\times {\mathbb R}) \simeq H^2 (W, {\mathbb Z}_G)\times H^2 (W, {\mathbb R}),
\end{gather*}
where ${\mathbb Z}_G:=\ker \{\exp\colon {\mathfrak g}\to G\}$ denotes the integral lattice of
the torus $G$.
\end{enumerate}
\end{Theorem}

The main dif\/f\/iculty in proving Theorem~\ref{thm:main} lies in establishing
part~(\ref{1}).
Results similar to part~(\ref{2}) hold in a somewhat greater generality
for completely integrable systems with elliptic singularities
(under a mild properness assumption) \cite{BM, Y, Z}:
once one knows that
there is {\it one} completely integrable system with the space of
tori $W$, the space of isomorphism classes of all such systems
is classif\/ied by the
second cohomology of~$W$ with coef\/f\/icients in an appropriate sheaf
(q.v.\ op.\ cit.).  The existence part for completely integrable
systems, called the realization problem by Zung~\cite{Z}, is much more
dif\/f\/icult.  For instance in~\cite{Z} the realization problem is only
addressed for 2-dimensional spaces of tori.  The solution to the
realization problem announced in~\cite{BM} and a similar solution in~\cite{Y} is dif\/f\/icult to apply in practice.
The solution, is, roughly, as follows.  Given an integral af\/f\/ine
manifold with corners $W$ one shows f\/irst that there is an open cover
$\{U_\alpha\}$ of~$W$ such that over each~$U_\alpha$ the realization
problem has a solution~$M_\alpha$.  Then, {\it if} there exists a
collection of isomorphisms $\varphi_{\alpha,\beta}\colon M_\beta|_{U_\alpha
  \cap U_\beta}\to M_\alpha|_{U_\alpha \cap U_\beta}$ satisfying the
appropriate cocycle condition, the realization problem has a solution
for~$W$.  Compare this with Theorem~\ref{thm:main}(\ref{1}) which
asserts that the realization problem for completely integrable torus
actions {\it always} has a solution.  We believe that the realization
problem for completely integrable systems with elliptic singularities
also always has a solution.  We will address this elsewhere.

Our proof of  part~(\ref{1}) of Theorem~\ref{thm:main} proceeds as follows.
We def\/ine symplectic toric $G$-bundles over~$\psi$:
these are
symplectic principal $G$-bundles over manifolds with corners
with orbital moment map~$\psi$.
They form a category, which we denote by ${\textsf{STB}}_\psi (W)$.
This category is always non-empty:
it contains the pullback $\psi^*(T^*G\to {\mathfrak g}^*)$.
We then construct a functor
\begin{gather} \label{collapse 0}
 {\textsf{c}} \colon \  {\textsf{STB}}_{\psi}(W)   \to
                   {\textsf{STM}}_{\psi}(W)
\end{gather}
from the category of symplectic toric $G$-bundles over $W$ to the
category ${\textsf{STM}}_{\psi}(W)$ of symplectic toric $G$-manifolds over~$W$.
The functor
${\textsf{c}}$ trades corners for f\/ixed points; it is a~version of a~symplectic cut~\cite{Lcuts}.  It follows, since
${\textsf{STB}}_\psi (W)$ is non-empty, that there always
exist symplectic toric $G$-manifolds
over a given unimodular local embedding $\psi\colon W \to {\mathfrak g}^*$
of a manifold with corners~$W$.

More is true.  We show that the functor ${\textsf{c}}$ is an
equivalence of categories. Hence, it induces a~bijection,
$\pi_0({\textsf{c}})$, between the isomorphism classes of objects of our
categories:
\begin{gather*}
\pi_0({\textsf{c}}) \colon \  \pi_0({\textsf{STB}}_\psi (W))   \to
                  \pi_0({\textsf{STM}}_\psi(W)).
\end{gather*}
The geometric meaning of the cohomology classes in $H^2 (W; {\mathbb Z}_G \times {\mathbb R})$
that classify symplectic toric $G$-manifolds over $W$ now becomes clear:
\emph{the elements of $H^2(W; {\mathbb Z}_G)$ classify principal $G$-bundles,
and the elements of $H^2(W; {\mathbb R})$
keep track of the ``horizontal part'' of the symplectic forms on these bundles.}

We note that, for compact symplectic toric $G$-manifolds,
the idea to obtain their classif\/ication by expressing these manifolds
as the symplectic cuts of symplectic toric $G$-manifolds with free~$G$ actions
is due to Eckhard Meinrenken (see~\cite[Chapter~7, Section~5]{Meinrenken}).

The paper is organized as follows.
In Section~\ref{sec:collapse}, after introducing our notation and
conventions, we construct the functor~\eqref{collapse 0}.  In
Section~\ref{sec:local uniqueness STB} we show that any two symplectic
toric~$G$ bundles over the same unimodular local embedding are locally
isomorphic (Lemma~\ref{loc-iso-stb}).
In Section~\ref{sec:equiv of cat} we
prove that the functor ${\textsf{c}}$ in~\eqref{collapse 0} is an equivalence of categories.  In
Section~\ref{sec:invariants}, we give the classif\/ication of symplectic
toric $G$-bundles over a f\/ixed unimodular local embedding $\psi\colon W \to
{\mathfrak g}^*$ in terms of two characteristic classes, the Chern class
$c_1$, which is in $H^2(W,{\mathbb Z}_G)$ and encodes the ``twistedness'' of
the $G$ bundle, and the horizontal class $c_{\text{hor}}$, which is in
$H^2(W,{\mathbb R})$ and encodes the ``horizontal part'' of the symplectic form
on the bundle.  We show that the map
\begin{gather*}
(c_1, c_{\text{hor}})\colon \  \pi_0 ({\textsf{STB}}_\psi(W)) \to H^2(W, {\mathbb Z}_G)\times H^2(W, {\mathbb R})
\end{gather*}
is a bijection.  Since the map $\pi_0({\textsf{c}}) \colon \pi_0({\textsf{STB}}_{\psi}(W))
\to \pi_0({\textsf{STM}}_{\psi}(W)) $ is a bijection, the composite
\begin{gather*}
\pi_0({\textsf{STM}}_\psi (W)) \xrightarrow{(c_1,c_{\text{hor}})\circ \pi_0(c)^{-1} }
H^2(W, {\mathbb Z}_G \times {\mathbb R})
\end{gather*}
is a bijection as well.
This classif\/ies (isomorphism classes of) symplectic toric
$G$-manifolds over $\psi\colon W \to {\mathfrak g}^*$.

Finally, in Section~\ref{sec:the_rest} we discuss those symplectic
toric manifolds that are determined by their moment map images.
In Proposition~\ref{Delzant over subset} we use Theorem~\ref{thm:main}
to derive Delzant's classif\/ication theorem and its generalization in the
case of symplectic toric $G$-manifolds that are not necessarily
compact but whose moment maps are proper as maps to convex subsets of~${\mathfrak g}^*$.
(In fact, already in~\cite{ka:appendix} it was noted that,
with the techniques of Condevaux--Dazord--Molino~\cite{CDM},
Delzant's proof should generalize to non-compact manifolds if the moment map
is proper as a map to a convex open subset of the dual of the Lie algebra.)
In Theorem~\ref{reduction of CN}, which was obtained in collaboration
with Chris Woodward,
we characterize those symplectic toric manifolds that
are symplectic quotients of the standard~${\mathbb C}^N$ by a subtorus of the
standard torus~${\mathbb T}^N$.
In Example~\ref{ex:2.15} we construct
a symplectic toric manifold that cannot be obtained by such a reduction.

The paper has two appendices.
Appendix~\ref{app:mfld-w-corner} contains background on manifolds with corners.
In Appendix~\ref{sec:background}, we recall the local normal form for
neighborhoods of torus orbits in symplectic toric manifolds,
and we use it to prove the following facts,
which are known but maybe hard to f\/ind in the literature:
\begin{enumerate}[1)]\itemsep=0pt
\item
orbit spaces of symplectic toric manifolds are
manifolds with corners;
\item
orbital moment maps of symplectic toric manifolds are
unimodular local embeddings; and
\item
any two symplectic toric manifolds over the same unimodular local embedding
are locally isomorphic.
\end{enumerate}

In the remainder of this section, following referees' suggestions,
we describe some relations of our work to existing literature
on integral af\/f\/ine structures and Lagrangian f\/ibrations.

\begin{Remark}[orbital moment maps]\label{add_rk:1}
  An {\em equivariant} moment map $\nu\colon N\to {\mathfrak h}^*$ for an action of
  a~Lie group $H$ on a symplectic manifold $N$ descends to a
  continuous map $\bar{\nu}\colon N/H\to {\mathfrak h}^*/H$ between orbit spaces.
  This map was introduced by Montaldi~\cite{Montaldi} under the name
  of {\em orbit momentum map} and was used to study stability and
  persistence of relative equilibria in Hamiltonian systems.  An
  analogue of this map in contact geometry was used by Lerman to
  classify contact toric manifolds~\cite{Ltoric}.

  The content of Proposition~\ref{prop:unimod-emb 1} is that for a
  symplectic toric manifold $(M,\omega, \mu)$ the orbit space $M/G$ is
  not just a topological space.  It has a natural structure of a
  $C^\infty$ manifold with corners and that the induced orbital moment
  map is $C^\infty$.

  Symplectic toric manifolds, in addition to being examples of
  symplectic manifolds with Hamiltonian torus actions, are also a
  particularly nice class of completely integrable systems with
  elliptic singularities.  Viewed this way $\bar{\mu}\colon M/G\to {\mathfrak g}^*$ is
  a developing map for an integral af\/f\/ine structure on the manifold
  with corners $M/G$ (see also Remark~\ref{add_rk:2} below).
\end{Remark}

\begin{Remark}[integral af\/f\/ine structures]
\label{add_rk:2}
An integral af\/f\/ine structure on a manifold with corners is usually
def\/ined in terms of an atlas of coordinate charts with integral af\/f\/ine
transition maps; see, for example,~\cite{Y}.  It is not hard to see
that such an atlas on a manifold~$W$ def\/ines  a~Lagrangian subbundle
${\mathcal L}$ of the cotangent bundle $T^*W\to W$ with two properties:
\begin{enumerate}[1)]\itemsep=0pt
\item the f\/iber ${\mathcal L}_w \subset T^*_wW$ is a lattice;
\item if $w\in W$ lies in a stratum of $W$ of
  codimension $k$ then there is a local frame $\{\alpha_1, \ldots,
  \alpha_n\}$ of $T^*W$ def\/ined near $w$ ($n=\dim W)$ so that the
  f\/irst $k$ 1-forms $\alpha_1,\ldots, \alpha_k$ annihilate the vectors
  tangent to the stratum.
\end{enumerate}
Conversely, any such Lagrangian subbundle def\/ines on $W$ an atlas of
coordinate charts with integral af\/f\/ine transition maps.

In general the  bundle ${\mathcal L}\to W$ may have no global
frame.  And even if it does have a~global frame $\{\alpha_1, \ldots,
\alpha_n\}$ the one forms $\alpha_j$ (which are necessarily closed)
need not be exact.  But if there is a global exact frame
$\{df_1,\ldots, df_n\}$ of ${\mathcal L}\to W$, then we have a smooth map
$f=(f_1,\ldots, f_n)\colon W\to {\mathbb R}^n$.  Such a~map~$f$ is a developing
map for the integral af\/f\/ine structure on~$W$.

Observe that a unimodular local embedding $\psi\colon W\to {\mathfrak g}^*$
def\/ines an integral af\/f\/ine structure on $W$ as follows. Since $\psi$
is a local embedding, the cotangent bundle $T^*W$ is the pullback by
$\psi$ of the cotangent bundle $T^*{\mathfrak g}^*$.  Consequently the standard
Lagrangian lattice ${\mathcal L}_{\rm can} ={\mathfrak g}^* \times {\mathbb Z}_G \subset {\mathfrak g}^*
\times {\mathfrak g}\simeq T^*{\mathfrak g}^*$ pulls back to a Lagrangian subbundle of
$T^*W$.  A choice of a basis of $\{e_1,\ldots, e_n\}$ of the integral
lattice ${\mathbb Z}_G$ def\/ines a map $f\colon W\to {\mathbb R}^n$.  It is given by
\begin{gather*}
f(w) = (\langle \psi, e_1\rangle, \ldots, \langle \psi, e_n\rangle).
\end{gather*}
The map $f$ is a developing map for $\psi^*{\mathcal L}_{\rm can}\to W$.
\end{Remark}

\subsection*{Symplectic toric manifolds and proper Lagrangian f\/ibrations}

Let $(M,\omega,\mu)$ be a~symplectic toric $G$-manifold
with a~$G$ quotient map $\pi \colon M \to W$.
Restricting to the interior ${\mathaccent23{W}}$ of~$W$ (as a manifold with corners),
we get a completely integrable system in the sense that was studied
by Duistermaat~\cite{Duis}, namely, a proper Lagrangian f\/ibration
with connected f\/ibers.
These were revisited and generalized by Dazord and Delzant~\cite{DD}.
For a detailed exposition see~\cite{LTB}.

\begin{Remark}[the integral af\/f\/ine structure and the monodromy]
As Duistermaat explains, a~proper Lagrangian f\/ibration
with connected f\/ibers $\pi \colon M \to B$
def\/ines an integral af\/f\/ine structure on the base $B$.  Each covector
$\beta \in T_b^*B$ determines a vector f\/ield $\xi_\beta$ along $\pi^{-1}(b)$
by the equation $\iota(\xi_\beta)\omega = \pi^* \beta$,
and the Lagrangian lattice sub-bundle is
\begin{gather*}
 {\mathcal L} =
   \{ \beta \, | \, \text{the f\/low of $\xi_\beta$ is $2\pi$ periodic} \}.
   \end{gather*}
Duistermaat's monodromy measures the non-triviality of the Lagrangian
lattice sub-bundle \mbox{${\mathcal L} \to B$}.
When it is trivial, the bundle of tori $T^*B/{\mathcal L} \to B$
becomes a trivial bundle with f\/iber, say,~$G$,
 $T^*B$ and ${\mathcal L}$ become trivial bundles
with f\/ibers ${\mathfrak g}^*$ and ${\mathbb Z}_G^*$,
and $\pi \colon M \to B$ becomes a $G$ principal bundle.
In this case, an orbital moment map is also
a developing map for the integral af\/f\/ine structure.  Having a moment map
in this context exactly means that the integral af\/f\/ine structure on~$B$
is globally developable.
\end{Remark}

\begin{Remark}[the characteristic classes] \label{rmrk:1.9}
Let $\pi \colon M \to B$ be a proper Lagrangian f\/ibration with connected f\/ibers.
The f\/ibers of the bundle of tori $T^*B/{\mathcal L}$ act freely and transitively
on the f\/ibers of $\pi \colon M \to B$.
Moreover, every point in $B$ has a neighborhood over which
$\pi \colon M \to B$ and $T^*B/{\mathcal L} \to B$ are isomorphic;
this is Duistermaat's formulation of the Arnold--Liouville theorem
on the local existence of action angle variables.
Globally, such f\/ibrations $\pi \colon M \to B$
are classif\/ied by the f\/irst cohomology group
\begin{gather*}
 H^1 \big( C^{\infty}_{\text{Lagr}} (\cdot,T^*B/{\mathcal L}) \big)
 \end{gather*}
of the sheaf of Lagrangian sections of $T^*B/{\mathcal L}$.

The short exact sequence of sheaves
\begin{gather*}
 0 \to C^{\infty}(\cdot,{\mathcal L}) \to C^{\infty}_{\text{Lagr}} (\cdot,T^*B)
  \to C^{\infty}_{\text{Lagr}}(\cdot,T^*B/{\mathcal L}) \to 0
\end{gather*}
gives an exact sequence
\begin{gather*}
 \cdots \to
   H^1\big(C^{\infty}_{\text{Lagr}}(\cdot,T^*B)\big) \to H^1\big(C^{\infty}_{\text{Lagr}}(\cdot,T^*B/{\mathcal L})\big)
   \to H^2(B,{\mathcal L}) \to \cdots .
\end{gather*}
Noting that Lagrangian sections of $T^*B$ are the same as closed one-forms,
and identifying the~$H^1$ of their sheaf with~$H^2(B,{\mathbb R})$,
we get an exact sequence
\begin{gather} \label{exact Duis}
 \cdots \to
 H^2(B,{\mathbb R}) \to H^1\big(C^{\infty}_{\text{Lagr}}(\cdot,T^*B/{\mathcal L})\big) \stackrel{c_1}{\to}
 H^2(B,{\mathcal L})
 \to \cdots.
\end{gather}

The second of these maps is Duistermaat's Chern class.  When the monodromy
is trivial, Duistermaat's Chern class is the Chern class
of $\pi \colon M \to B$ as a principle $G$ bundle.
If $\pi \colon M \to W$ is the $G$-quotient map of a symplectic toric
$G$-manifold, then Duistermaat's Chern class for $M|_{{\mathaccent23{W}}}$
coincides with ours under the identif\/ication
$H^2(W;{\mathbb Z}_G) \stackrel{\cong}{\to} H^2({\mathaccent23{W}};{\mathbb Z}_G)$.

If the monodromy and Chern class both vanish, Duistermaat def\/ines a class
in $H^2(B,{\mathbb R})$, which is often called the \textit{Lagrangian class};
it is the cohomology class of the pullback of $\omega$ by a~global smooth
section.  If $\pi \colon M \to B$ is the $G$ quotient map of a symplectic toric
$G$-manifold, and if additionally the Chern class vanishes,
then Duistermaat's Lagrangian class for $M|_{{\mathaccent23{W}}}$ coincides
with our horizontal class under the identif\/ication
$H^2(W;{\mathbb R}) \stackrel{\cong}{\to} H^2({\mathaccent23{W}};{\mathbb R})$.
\end{Remark}

\begin{Remark}\label{rmrk:1.11}
If $\pi \colon M \to W$ is the $G$-quotient map for a symplectic toric
$G$-manifold and $B = {\mathaccent23{W}}$, our characteristic class gives a splitting
\begin{gather*}
 H^1\big(C^{\infty}_{\text{Lagr}}(B,T^*B/{\mathcal L})\big) \cong H^2\big(B;{\mathbb Z}_G^*\big) \oplus H^2(B;{\mathbb R})
\end{gather*}
that is consistent with~\eqref{exact Duis}.
Moreover, our construction
provides a geometric meaning to the Lagrangian class in~$H^2(B,{\mathbb R})$.

In this case
\begin{itemize}\itemsep=0pt
\item
every element of $H^2(W;{\mathbb Z}_G)$ gives rise to a symplectic toric $G$-manifold,
and
\item
distinct elements of $H^2(W;{\mathbb R})$ represent non-isomorphic symplectic toric
$G$-manifolds.
\end{itemize}
Both of these facts are not necessarily true in the more general situation
that is addressed by Duistermaat and Dazord--Delzant.
\end{Remark}

\section{A functor from symplectic toric bundles\\ to symplectic toric manifolds}
\label{sec:collapse}

The purpose of this section is to construct a  functor
 \begin{gather*}
 {\textsf{c}} \colon \  {\textsf{STB}}_{\psi}(W)   \to
                   {\textsf{STM}}_{\psi}(W)
\end{gather*}
from the category of symplectic toric $G$-bundles to the category of
symplectic toric $G$-manifolds over a given unimodular local embedding
$\psi\colon W\to {\mathfrak g}^*$ of a manifold with corners $W$. Once this functor
is constructed, we deduce Theorem~\ref{thm:main}(\ref{1}) almost
immediately.  In Section~\ref{sec:equiv of cat} we  prove that~${\textsf{c}}$ is an equivalence of categories.  We start by establishing
our notation and recording a few necessary
def\/initions.

\textbf{Notation and conventions.}
A \emph{torus} is a compact connected abelian Lie group.  A torus
of dimension $n$ is isomorphic, as a Lie group, to $(S^1)^n$
and to ${\mathbb R}^n/{\mathbb Z}^n$.  We denote the Lie algebra of a torus $G$ by ${\mathfrak g}$,
the dual of the Lie algebra, $\text{Hom}({\mathfrak g},{\mathbb R})$, by ${\mathfrak g}^*$, and
the integral lattice, $\ker (\exp \colon {\mathfrak g} \to G)$, by ${\mathbb Z}_G$.
The \emph{weight lattice} of $G$ is the lattice dual to ${\mathbb Z}_G$;
we denote it by ${\mathbb Z}_G^*$.
When a torus $G$ acts on a manifold $M$, we denote the action
of an element $g \in G$ by $m \mapsto g \cdot m$ and the vector f\/ield
induced by a Lie algebra element $\xi \in {\mathfrak g}$ by $\xi_M$; by
def\/inition
\begin{gather*}
\xi_M (m) = \left. \frac{d}{dt}\right|_{t=0} (\exp
(t\xi)\cdot m).
\end{gather*}
We write the canonical pairing between~${\mathfrak g}^*$ and~${\mathfrak g}$
as $\langle \cdot, \cdot\rangle$.  Our sign convention for a moment
map $\mu \colon M \to {\mathfrak g}^*$ for a Hamiltonian action of a torus~$G$
on a symplectic manifold~$(M, \omega)$ is that it satisf\/ies
\begin{gather} \label{eq:moment_map}
d \langle \mu, \xi \rangle = - \omega (\xi_M, \cdot) \qquad \textrm{for all} \ \ \xi \in {\mathfrak g}.
\end{gather}
For us a \emph{symplectic toric $G$-manifold} is a triple $(M,\omega, \mu)$
where $M$ is a manifold, $\omega$ is a symplectic form and
$\mu\colon M\to{\mathfrak g}^*$ is a moment map for an ef\/fective Hamiltonian action of
a torus $G$ with $\dim M = 2 \dim G$.

\begin{Definition} \label{def:unimod_cone}
A \emph{unimodular cone} in the dual ${\mathfrak g}^*$ of the Lie algebra
of a torus $G$ is a subset~$C$ of~${\mathfrak g}^*$ of the form
\begin{gather*} 
C= \big\{ \eta \in {\mathfrak g}^* \,|\, \langle \eta - \epsilon , v_i \rangle \geq 0
\; \text{for all}\; 1 \leq i \leq k\big\},
\end{gather*}
where $\epsilon$ is a point in ${\mathfrak g}^*$ and $\{v_1, \ldots, v_k\}$ is a
basis of the integral lattice of a subtorus of $G$.  We record the
dependence of the cone $C$ on the data $\{v_1,\ldots,v_k\}$ and $\epsilon$ by
writing
\begin{gather*}
C = C_{\{v_1,\ldots,v_k\}, \epsilon}.
\end{gather*}
\end{Definition}

\begin{Remark}
The set $C={\mathfrak g}^*$ is a unimodular cone def\/ined by the empty basis
of the integral lattice~$\{ 0 \}$ of the trivial subtorus~$\{1\}$ of~$G$.
\end{Remark}

\begin{Remark}
A unimodular cone is a manifold with corners.
Moreover, it is a manifold with faces (q.v.\ Def\/inition~\ref{def:faces}).

For a unimodular cone $C= C_{\{v_1,\ldots,v_k\}, \epsilon}$ the
facets are the sets
\begin{gather*}
F_i= \{ \eta \in C\,|\, \langle \eta -\epsilon , v_i\rangle =0\},
        \qquad 1\leq i\leq k.
\end{gather*}
The vector $v_i$ in the formula above is the inward
pointing primitive normal to the facet~$F_i$. (Recall that a
vector $v$ in the lattice ${\mathbb Z}_G$ is {\em primitive}
if for any $u\in {\mathbb Z}_G$ the equation $v=nu$ for $n\in {\mathbb Z}$
implies that $n= \pm 1$.)
\end{Remark}

\begin{Lemma}\label{lem:prim-norm}
  The primitive inward pointing normal $v_i$ to a facet $F_i$ of a~unimodular cone $C_{\{v_1,\ldots,v_k\}, \epsilon}$ is uniquely
  determined by any open neighborhood ${\mathcal O}$ of a point~$x$ of~$F_i$
  in~$C$.
\end{Lemma}

\begin{proof}
The af\/f\/ine hyperplane spanned by $F_i$ is uniquely
determined by the intersection \mbox{${\mathcal O}\cap F_i$}.  Up to sign, such a
hyperplane has a unique primitive normal.  The sign of the normal is
determined by requiring that at the point $x$ the normal points into
${\mathcal O}$.
\end{proof}

\begin{Definition}[unimodular local embedding (u.l.e.)]
\label{def:unimodular emb}
  Let $W$ be a manifold with corners and~${\mathfrak g}^*$ the dual of
  the Lie algebra of a torus.  A smooth map $\psi \colon W \to {\mathfrak g}^*$
  is a~\emph{unimodular local embedding} (a~{\em u.l.e.}) if for each
  point~$w$ in~$W$ there exists an open neighborhood ${\mathcal T} \subset W$ of
  the point and a unimodular cone $C \subset {\mathfrak g}^*$ such that
  $\psi({\mathcal T})$ is contained in $C$ and $\psi|_{{\mathcal T}} \colon {\mathcal T}
  \to C$ is an {\em open embedding}. That is, $\psi({\mathcal T})$ is open in $C$
  and $\psi|_{{\mathcal T}}\colon  {\mathcal T}\to \psi({\mathcal T})$ is a dif\/feomorphism.
\end{Definition}

\begin{Remark}
In Def\/inition~\ref{def:unimodular emb},
the cone $C$ is not uniquely determined by the point $w$;
for instance it can have facets that do not pass through $\psi(w)$.
For example, let $G = (S^1)^2$, let $\psi \colon W \to {\mathfrak g}^* = {\mathbb R}^2$
be the inclusion map of the positive quadrant, and let $w = (1,0)$.
If the neighborhood ${\mathcal T}$ of $w$ meets the non-negative $y$ axis,
then the cone $C$ must be the positive quadrant too.  Otherwise,
the natural choice for $C$ is the closed upper half plane,
but for suitable choices of ${\mathcal T}$ the cone $C$ can also be
the intersection of the closed upper half plane
with a half plane of the form $\{ x+ny \geq c \}$ for $n \in {\mathbb Z}$ and $c<1$
or of the form $\{ x+ny \leq c \}$ for $n \in {\mathbb Z}$ and $c>1$.
\end{Remark}

\begin{Remark}
  Proposition~\ref{prop:unimod-emb 1} shows that the orbital moment
  map of a symplectic toric manifold is a unimodular local embedding.
\end{Remark}

\begin{Example} \label{ex:mmap-not-emb}
  It is easy to construct examples where the orbital moment map is not
  an embedding.  Consider, for instance, a 2-dimensional torus $G$.
  Removing the origin from the dual of its Lie algebra ${\mathfrak g}^*$ gives
  us a space
   that is
  homotopy equivalent to a circle.  Thus the f\/ibers of the
  universal covering map $p \colon W\to {\mathfrak g}^*\setminus \{0\}$ have
  countably many points.  The pullback $p^*(T^*G)$ along $p$ of the
  principal $G$-bundle $\mu \colon T^*G \to {\mathfrak g}^*$ is a symplectic toric
  $G$-manifold with orbit space~$W$ and orbital moment map~$p$, which is
  certainly not an embedding.

Similarly, let $S^2$ be the unit sphere in ${\mathbb R}^3$ with the standard
area form, and equip~$S^2 \times S^2$ with the standard toric action
of $(S^1)^2$ with moment map $\mu((x_1,x_2,x_3),(y_1,y_2,y_3)) = (x_3,y_3)$.
Its image is the square $I^2 = [-1,1] \times [-1,1]$.
Remove the origin, and let
$p \colon W \to I^2 \setminus \{ 0 \}$ be the universal covering.
Then the f\/iber product $W \times_{I^2 \setminus \{ 0 \} } (S^2 \times S^2)$
is a symplectic toric manifold; it is a ${\mathbb Z}$-fold covering
of $(S^2 \times S^2) \setminus ( \text{the equator} \times \text{the equator})$.
As in the previous example, the orbital moment map is not an embedding.
Unlike the previous example, the torus action  here is not free.
\end{Example}

\begin{Definition}
  Let $W$ be a manifold with corners and $\psi\colon W\to{\mathfrak g}^*$ a unimodular
  local embedding.
  A {\em symplectic toric manifold over $\psi\colon W\to{\mathfrak g}^*$}
  is a symplectic toric $G$-manifold $(M,\omega,\mu)$,
  equipped with a quotient map $\pi\colon M\to W$
  for the action of $G$ on $M$ (q.v.\ Def\/inition~\ref{quotient}),
  such that
\begin{gather*}
\mu = \psi \circ \pi.
\end{gather*}
\end{Definition}

\begin{Remark}
  Since the moment map $\mu\colon M\to {\mathfrak g}^*$ together with the symplectic
  form $\omega$ encodes the action of the group $G$ on $M$ and since
  the quotient map $\pi\colon M\to W$ together with $\psi\colon W\to{\mathfrak g}^*$ encode
  $\mu$, we may regard a symplectic toric $G$-manifold over $\psi\colon W\to
  {\mathfrak g}^*$ as a triple $(M,\omega, \pi\colon M\to W)$.
\end{Remark}

We now f\/ix a u.l.e.\
\begin{gather*}
\psi \colon \  W \to {\mathfrak g}^*
\end{gather*}
of a manifold with corners $W$ into the dual of the Lie algebra
of a torus $G$,
and proceed to def\/ine the category ${\textsf{STM}}_\psi(W)$
of \emph{symplectic toric $G$-manifolds} over $W\stackrel{\psi}{\to} {\mathfrak g}^*$
and the category ${\textsf{STB}}_\psi(W)$
of \emph{symplectic toric $G$-bundles} over $W\stackrel{\psi}{\to} {\mathfrak g}^*$.

\begin{Definition}[the category $\mathbf{{\textsf{STM}}_\psi (W)}$ of symplectic toric
  $G$-manifolds over $\mathbf{\psi\colon W\to {\mathfrak g}^*}$] \label{def:STM}
  We
  def\/ine an \emph{object} of the category ${\textsf{STM}}_\psi(W)$ to be a~symplectic toric $G$-manifold  \mbox{$(M,\omega, \pi\colon M{\to} W)$} over~$W$.
  A \emph{morphism} $\varphi$ from $(M, \omega,  \pi)$
  to $(M', \omega', \pi') $
  is a $G$-equivariant symplectomorphism $\varphi \colon M \to M'$
  such that $\pi' \circ \varphi = \pi$.
\end{Definition}

\begin{Notation} \label{notation for STM}\sloppy
Informally, we may sometimes write  $M$ for an object of
${\textsf{STM}}_\psi (W)$ and \mbox{$\varphi \colon M \to M'$} for a
morphism between two such objects.
Also, we may write $W$ as shorthand for
\mbox{$\psi \colon W \to {\mathfrak g}^*$}.
\end{Notation}

\begin{Definition}[the category $\mathbf{{\textsf{STB}}_\psi (W)}$ of symplectic toric
  $G$-bundles over $\mathbf{\psi\colon W\to {\mathfrak g}^*}$] \label{def:STB}
  An  \emph{object} of the category ${\textsf{STB}}_\psi (W)$ is a principal
  $G$-bundle $\pi\colon P\to W$ over a manifold with corners (cf.\
  Def\/inition~\ref{def:pric_G-bundle_corners}) together with a~$G$-invariant symplectic form $\omega$ so that $\mu := \psi \circ
  \pi$ is a moment map for the action of~$G$ on~$(P,\omega)$.
  We call the triple $(P,\omega, \pi\colon  P \to W)$
  a {\em symplectic toric $G$-bundle over the u.l.e.\ $\psi\colon  W\to {\mathfrak g}^*$},
  or a symplectic toric $G$-bundle over~$W$ for short.
  A~\emph{morphism}~$\varphi$ from $(P, \omega, \pi)$ to
  $(P', \omega', \pi')$ is a~$G$-equivariant
  symplectomorphism~$\varphi \colon P \to P'$ with $\pi' \circ
  \varphi = \pi$.
\end{Definition}

\begin{Remark}
 The categories ${\textsf{STM}}_\psi(W)$ and ${\textsf{STB}}_\psi(W)$
are groupoids, that is, all of their morphisms are invertible.

  If $\psi \colon W\to {\mathfrak g}^*$ is a u.l.e., $(M,\omega, \pi)$ is a symplectic
  toric $G$ manifold over $W$ and $U\subset W$ is open, then the restriction
  $\psi|_U\colon U \to {\mathfrak g}^*$ is also a u.l.e., and
\begin{gather*}
\big(M|_U:= \pi^{-1} (U),\, \omega|_{M|_U},\, \pi|_{M|_U}\big)
\end{gather*}
is a symplectic toric $G$-manifold over $U$.
The restriction map extends to a functor
\begin{gather*}
|^W_U\colon \  {\textsf{STM}}_\psi (W) \to {\textsf{STM}}_\psi(U).
\end{gather*}
Given an open subset $V$ of $U$ we get the restriction $|_V^U\colon
{\textsf{STM}}_\psi(U)\to {\textsf{STM}}_\psi(V)$.  The three restriction functors are
compatible:
\begin{gather*}
|^W_V = |^U_V \circ |^W_U.
\end{gather*}
In other words, the assignment
\begin{gather*}
U\mapsto {\textsf{STM}}_\psi (U)
\end{gather*}
is a (strict) presheaf of groupoids.

A reader familiar with stacks will have little trouble checking that
the presheaf ${\textsf{STM}}_\psi (\cdot)$ satisf\/ies the descent condition with
respect to any open cover of~$W$ and that thus ${\textsf{STM}}_\psi (\cdot )$ is
a~stack on the site ${\rm Open}(W)$ of open subsets of~$W$ with the cover
topology.  The stack~${\textsf{STM}}_\psi $ is not a~geometric stack.

Similarly, a symplectic toric bundle over a
manifold with corners~$W$ restricts to a symplectic toric bundle over
an open subset of~$W$.  These restrictions def\/ine a presheaf of
groupoids ${\textsf{STB}}_\psi (\cdot )$.  A~reader familiar with stacks can
check that ${\textsf{STB}}_\psi (\cdot )$ is also a stack; see also
Lemma~\ref{lem:descent_for_stb} below.
\end{Remark}

\begin{Remark}\label{interior}
  If $\psi \colon W\to {\mathfrak g}^*$ is a u.l.e.\ and $W$
is a manifold without corners (i.e., a manifold) then
\begin{gather*}
{\textsf{STM}}_\psi(W) = {\textsf{STB}}_\psi (W).
\end{gather*}
If $W$ is an arbitrary manifold with corners, then its interior ${\mathaccent23{W}}$
(q.v. Def\/inition~\ref{strata}) is a~manifold, and so
\begin{gather*}
{\textsf{STM}}_\psi({\mathaccent23{W}}) = {\textsf{STB}}_\psi ({\mathaccent23{W}}).
\end{gather*}
\end{Remark}

\subsection*{The functor $\boldsymbol{{\textsf{c}} \colon{\textsf{STB}}_{\psi}(W)  \to   {\textsf{STM}}_{\psi}(W)}$}

Next we outline the construction of the functor ${\textsf{c}} \colon
  {\textsf{STB}}_{\psi}(W)   \to   {\textsf{STM}}_{\psi}(W) $ from the category of
symplectic toric $G$-bundles to the category of symplectic toric
$G$-manifolds over a u.l.e.~$\psi$.

{\bf Step 1: characteristic subtori.} We show that $\psi$
attaches to each point $w\in W$ a subto\-rus~$K_w$ of~$G$ together with
a choice of a basis
$\{v_1^{(w)},\ldots, v_k^{(w)}\}$ of its integral lattice ${\mathbb Z}_{K_w}$.

A basis of the integral lattice ${\mathbb Z}_K$ of a torus $K$
def\/ines a linear symplectic representation
$\rho\colon K\to \operatorname{Sp}(V,\omega_V)$, which we may regard as a symplectic toric
$K$-manifold $(V,\omega_V, \mu_V)$ (here $\mu_V\colon V\to{\mathfrak k}^*$ is the
associated moment map with $\mu_V(0)=0$).  Thus for each point
$w\in W$ we also have a symplectic toric $K_w$-manifold $(V_w,
\omega_w, \mu_w)$.

  {\bf Step 2: a topological version $\boldsymbol{{\textsf{c}_{\rm top}}}$ of the
  functor $\boldsymbol{{\textsf{c}}}$.} The collection of the subtori $\{K_w\}_{w\in
  W}$ def\/ines for each principal $G$-bundle $\pi\colon P\to W$ an
equivalence relation $\sim$ in a functorial manner.  We show that
\begin{enumerate}\itemsep=0pt
\item Each quotient ${\textsf{c}_{\rm top}} (P):= P/{\sim}$ is a topological $G$-space
with orbit space $W$ and the action of $G$ on ${\textsf{c}_{\rm top}}(P)$
is free over the interior ${\mathaccent23{W}}$.
(Here, ${\textsf{c}_{\rm top}}$ stands for ``topological cut''.)
\item For every map $\varphi\colon P\to P'$ of principal $G$-bundles over
  $W$ we naturally get a $G$-equivariant homeomorphism
  ${\textsf{c}_{\rm top}}(\varphi)\colon {\textsf{c}_{\rm top}}(P)\to {\textsf{c}_{\rm top}}(P')$.
\item These data def\/ine a functor
\begin{gather*}
{\textsf{c}_{\rm top}}\colon \  {\textsf{STB}}_\psi (W)\to \text{topological $G$-spaces over $W$}.
\end{gather*}
\item Moreover, ${\textsf{c}_{\rm top}}$  is a map of presheaves of groupoids.
      In particular, for every open subset~$U$ of~$W$,
\begin{gather*}
{\textsf{c}_{\rm top}}(P|_U)= {\textsf{c}_{\rm top}}(P)|_U.
\end{gather*}
\end{enumerate}

{\bf Step 3: the actual construction of $\boldsymbol{{\textsf{c}}}$.} We
show that for every point $w\in W$ there is an open neighborhood
$U_w$ so that for every symplectic toric $G$-bundle $(P,\omega, \pi\colon P\to
W)$ the symplectic quotient
\begin{gather*}
\operatorname{cut}(P|_{U_w}):= (P|_{U_w} \times V_w)/\!/_0  K_w
\end{gather*}
is a symplectic toric $G$-manifold over $U_w$.

As in Step 2 the mapping $\operatorname{cut} ({\cdot}|_{U_w})$ (i.e.,
the restriction to~$U_w$ followed by $\operatorname{cut}$) from symplectic toric
$G$-bundles over $W$ to symplectic toric manifolds over $U_w$ extends
to a functor. In particular for every map $\varphi\colon P\to P'$ of
symplectic toric $G$-bundles over $W$ we have a map $\operatorname{cut}
(\varphi|_{U_w})\colon  \operatorname{cut}(P|_{U_w})\to \operatorname{cut}(P'|_{U_w})$ of
symplectic toric $G$-manifolds over $U_w$.

At the same time, for each symplectic toric $G$-bundle $P\to W$
we construct a collection
\begin{gather*}
\big\{\alpha_w^P\colon \; {\textsf{c}_{\rm top}}(P|_{U_w}) \to \operatorname{cut}(P|_{U_w})\big\}_{w\in W}
\end{gather*}
of equivariant homeomorphisms that have the following two
compatibility properties:
\begin{enumerate}\itemsep=0pt
\item  For a f\/ixed bundle $P\in {\textsf{STB}}_\psi(W)$ and any two points $w_1$, $w_2$
  the map
\begin{gather*}
\big(\alpha^P_{w_2 }\big)\circ \big(\alpha^P_{w_1}\big)^{-1} \colon \
\operatorname{cut}(P|_{U_{w_1}})|_{U_{w_1}\cap U_{w_2}}\to
\operatorname{cut}(P|_{U_{w_2}})|_{U_{w_1}\cap  U_{w_2}}
\end{gather*}
is a map of  symplectic  toric $G$-manifolds over $U_{w_1}\cap  U_{w_2}$.

\item For a point $w\in W$ and a map $\varphi\colon P_1\to P_2$
of symplectic toric bundles over $W$ the diagram
\begin{gather}\label{eq:alpha-cut}
\begin{split} &
\xymatrix{
 {\textsf{c}_{\rm top}}(P_1)|_{U_w} \ar[rr]^{\alpha^{P_1}_w} \ar[d]^{{\textsf{c}_{\rm top}}(\varphi)|_{U_w}}
 && \operatorname{cut}(P_1|_{U_w})\ar[d]^{\operatorname{cut}{(\varphi|_{U_w})}}\\
 {\textsf{c}_{\rm top}}(P_2)|_{U_w} \ar[rr]^{\alpha^{P_2}_w}
 && \operatorname{cut}(P_2|_{U_w})
}\end{split}
\end{gather}
commutes.
\end{enumerate}
The f\/irst property tells us that the family $\{\alpha_w^P\}_{w\in W}$
of homeomorphisms def\/ines on ${\textsf{c}_{\rm top}}(P)$ the structure of a symplectic
toric $G$-manifold over $\psi\colon W\to {\mathfrak g}^*$. We denote this manifold,
which is an object of ${\textsf{STM}}_\psi(W)$, by~$\textsf{c}(P)$.  The second property
tells us that ${\textsf{c}_{\rm top}}(\varphi)$ def\/ines a map $\textsf{c}(\varphi)\colon \textsf{c} (P_1)\to
\textsf{c}(P_2)$ of symplectic toric $G$-manifolds over~$\psi$.  This gives
rise to the
desired functor~$\textsf{c}$.

We now proceed to f\/ill in the details of the construction.
\begin{proof}[Details of Step~1] We start by proving

\begin{Lemma} \label{lemma:step1.1} Given a unimodular local embedding
  $\psi\colon W\to {\mathfrak g}^*$ and a point $w \in W$ there exists a unique subtorus $K_w$
  of $G$ and a unique basis $\{v_1^{(w)},\ldots, v_k^{(w)}\}$ of its integral
  lattice ${\mathbb Z}_{K_w}$ such that the following holds.
There exists an open neighborhood ${\mathcal U}_w$ of $w$ in $W$ so that
\begin{gather*}
\psi|_{{\mathcal U}_w}\colon \ {\mathcal U}_w \to C_w:= \big\{
\eta \in {\mathfrak g}^* \,|\,
\big\langle \eta -\psi (w),  v_j^{(w)}\big\rangle\geq 0 \text{ for } 1\leq j\leq k\big\}
\end{gather*}
is an open embedding of manifolds with corners.
\end{Lemma}

\begin{proof}
By def\/inition of a u.l.e., there exists an open neighborhood ${\mathcal T}
\subset W$ of $w$ and a~uni\-mo\-du\-lar cone $C = C_{\{u_1,\ldots,u_n\},
  \epsilon} \subset {\mathfrak g}^*$ such that $\psi({\mathcal T})$ is contained in $C$ and
$\psi|_{{\mathcal T}} \colon {\mathcal T} \to C$ is an open embedding of manifolds
with corners.  Since $\psi|_{{\mathcal T}}$ is an open embedding it maps the
interior of~${\mathcal T}$ to an open subset of the interior of~$C$.  We may
assume that~${\mathcal T}$ is a neighborhood with faces.  Then the stratum~$S$ of~${\mathcal T}$ containing~$w$ lies in exactly~$k$
facets ${\mathcal F}_1,\ldots, {\mathcal F}_k$ of~${\mathcal T}$, where~$k$ is the
codimension of~$S$.  For each~$j$ the image~$\psi({\mathcal F}_j)$ is an
open subset of a unique facet~$F_{i(j)}$ of~$C$ and~$\psi({\mathcal T})$ is
an open neighborhood of $\psi({\mathcal F}_j)$ in~$C$.
By Lemma~\ref{lem:prim-norm} the pair $(\psi({\mathcal T}), \psi({\mathcal F}_j))$ uniquely
determines the primitive inward pointing normal~$u_{i(j)}$ of the facet~$F_{i(j)}$ of~$C$.  Since $\{u_1,\ldots,u_n\}$ is a basis of an integral
lattice of a subtorus of~$G$, its subset $\{u_{i(j)}\}_{j=1}^k$ is
also a basis of an integral lattice of a possibly smaller subtorus
$K_w$ of $G$. We set $v_j^{(w)}:= u_{i(j)} $, $1\leq j\leq k$.
We  note that
\begin{gather*}
K_w = \exp\big(\operatorname{span}_{\mathbb R} \big\{v_1^{(w)},\ldots, v_k^{(w)}\big\}\big).
\end{gather*}
To obtain the neighborhood ${\mathcal U}_w$ we delete from the manifold with
faces~${\mathcal T}$ all the faces that {\em do not} contain~$w$.
\end{proof}

\begin{Remark}
  The basis $\{v_1^{(w)},\ldots, v_k^{(w)}\}$
  and the corresponding torus $K_w$ do not depend on our choice of the
  cone~$C$:  by construction $v_j^{(w)}$ is the primitive
  normal to the af\/f\/ine hyperplane spanned by $\psi ({\mathcal F}_j)$ that
  points into~$\psi({\mathcal T})$.  In fact the only way we use the
  existence of the unimodular cone $C$ is to insure that the set
  $\{v_1^{(w)},\ldots, v_k^{(w)}\}$ of normals to the facets of
  $\psi({\mathcal T})$ forms a basis of an integral lattice of a subtorus of
  the torus~$G$.

  Similarly, the basis $\{v_i^{(w)}\}_{i=1}^k$ does not depend on the choice of
  ${\mathcal T}$ either.
\end{Remark}

\begin{Remark}
  For each stratum of $W$ the function $w\mapsto K_w$ is
  locally constant, hence constant.  Consequently the subtorus~$K_w$
  depends only on the stratum of $W$ containing the point~$w$ and
  not on the point~$w$ itself.  Similarly the basis
  $\{v_1^{(w)},\ldots, v_k^{(w)}\}$ depends only on the stratum of~$W$
  containing $w$.
\end{Remark}

\begin{Remark}\label{rmrk:2.21}
For $w'\in {\mathcal U}_w$ we can read of\/f the group $K_{w'}$ from the face structure
of ${\mathcal U}_w$ and the set $\{v_1^{(w)},\ldots, v_k^{(w)}\}$. Namely
\begin{gather*}
K_{w'} = \exp \big(\operatorname{span}_{\mathbb R} \big\{v_i^{(w)} \, | \,
   \big\langle \psi(w')-\psi(w), v_i^{(w)}\big\rangle  = 0\big\}\big).
\end{gather*}
We also note that the subset
\begin{gather*}
\big\{v_i^{(w)} \, | \,  \big\langle \psi(w')-\psi(w), v_i^{(w)} \big\rangle = 0\big\}
\end{gather*}
of $\{v_1^{(w)},\ldots, v_k^{(w)}\}$ forms a basis of the integral
lattice of $K_{w'}$.
\end{Remark}

\begin{Lemma}\label{lem:mflds-w-faces}
A manifold with corners~$W$ that admits a u.l.e.\ $\psi \colon W \to {\mathfrak g}^*$
is a manifold with faces
$($q.v.~{\rm \cite{Janich}} and Definition~{\rm \ref{def:faces}} below$)$.
In particular, for any symplectic toric $G$-manifold $(M,\omega,\mu)$,
the quotient $M/G$ is a manifold with faces.
\end{Lemma}

\begin{proof}
  The map $\psi\colon W\to {\mathfrak g}^*$  sends a
  neighborhood of a point in a codimension~1 stratum $S$ of $W$ to a~relatively open subset of an af\/f\/ine hyperplane $H\subset {\mathfrak g}^*$
  whose normal lies the integral lattice ${\mathbb Z}_G$ of $G$.  Consequently
  $\psi$ sends all of $S$ to $H$ and $\psi|_S\colon S\to H$ is a~local dif\/feomorphism.  The lemma follows from this observation.
\end{proof}

\begin{Remark}
  It follows from the proof of Lemma~\ref{lem:mflds-w-faces} that the
  map $\psi\colon W\to {\mathfrak g}^*$ attaches to every (connected) codimesion~1
  stratum~$S$ of~$W$ a primitive vector $\lambda(S)\in {\mathbb Z}_G$
  (namely, the corresponding primitive inward normal).
  The function $S\mapsto \lambda(S)$ is the
  analogue of the characteristic function of Davis and Januszkiewicz~\cite{DJ} and of the characteristic bundle of Yoshida~\cite{Y}.
\end{Remark}

Recall that any symplectic representation of a torus is complex hence
has well-def\/ined weights.  These weights do not depend on a choice of
an invariant complex structure compatible with the symplectic form
since the space of such structures is path connected.

\begin{Lemma}\label{lemma:step1.2}
  Let $\rho_i\colon K\to \operatorname{Sp}(V_i,\omega_i)$, $i=1,2$, be two symplectic
  representations of a torus $K$ with the same set of weights.
  Then there exists a symplectic linear isomorphism of representations
  $\varphi\colon  (V_1,\omega_1) \to (V_2,\omega_2)$.
\end{Lemma}

\begin{proof}
  Choose $K$-invariant compatible complex structures on~$V_1$ and~$V_2$.
  As complex $K$ representations, each of~$V_1$ and~$V_2$
  decomposes into one-dimensional complex representations.
  Because the weights are the same, it is enough to consider the case
  that~$V_1$ and~$V_2$ are the same complex vector space
  and its complex dimension is one.
  In this case, because~$\omega_1$ and $\omega_2$ are both compatible
  with the complex structure, one must be a positive multiple of the other:
  $\omega_2 = \lambda^2 \omega_1$ for some scalar $\lambda >0$.
  We may then take $\varphi(v) := \lambda v$.
\end{proof}

Lemmas~\ref{lemma:step1.1} and \ref{lemma:step1.2} imply that to any
point $w$ of a manifold with corners $W$ a u.l.e.\ $\psi\colon W\to {\mathfrak g}^*$
unambiguously attaches a symplectic toric $K_w$-manifold $(V_w,
\omega_w, \mu_w)$: the weights of the representation $V_w$ is the
basis $\{v_j^*\}$ of the weight lattice ${\mathbb Z}_{K_w}^*$ dual to the basis
$\{ v_j^{(w)}\}$.  If $V'_w$ is another symplectic representation of
$K_w$ with the same set of weights as~$V_w$ then the symplectic toric
$K_w$-manifolds $(V_w, \omega_w, \mu_w)$ and $(V'_w, \omega_w',
\mu'_w)$ are linearly isomorphic as symplectic toric manifolds.
\end{proof}

\begin{proof}[Details of Step 2]\sloppy
  Given a principal $G$-bundle $\pi\colon P\to W$ we def\/ine $\sim$ to be the
  smallest equivalence relation on $P$ such that $p\sim p'$ whenever
  $\pi(p) = \pi(p')$ and $p$, $p'$ lie on the same~$K_{\pi(p)}$ orbit.
We give the set $P/_{\sim}$ the quotient topology.  Since the
  action of~$K_w$ on the f\/iber of~$P$ above $w$ commutes with the
  action of~$G$, the topological space
\begin{gather*}
{\textsf{c}_{\rm top}}(P) := P/_{\sim}
\end{gather*}
is naturally a $G$-space.  For
  the same reason $\pi\colon P\to W$ descends to a quotient map
  \mbox{$\bar{\pi}\colon {\textsf{c}_{\rm top}}(P) \to  W$}.  Since for the points $w$ in the interior
  of~$W$ the groups~$K_w$ are trivial, the action of~$G$ on~${\textsf{c}_{\rm top}}
  (P)|_{\mathaccent23{W}}$ is free.

  If $\varphi\colon P\to P'$ is a map of principal $G$-bundles over $W$,
  then it maps f\/ibers to f\/ibers and $K_w$-orbits to $K_w$ orbits
  thereby inducing ${\textsf{c}_{\rm top}}(\varphi)\colon {\textsf{c}_{\rm top}}(P)\to {\textsf{c}_{\rm top}}(P')$.
Explicitly ${\textsf{c}_{\rm top}}(\varphi)$ is given~by
\begin{gather*}
{\textsf{c}_{\rm top}}(\varphi)([p]) = [\varphi(p)].
\end{gather*}
Here, as before $[p]\in P/_{\sim} = {\textsf{c}_{\rm top}}(P)$ denotes the equivalence
class of $p\in P$ and $ [\varphi(p)]$ denotes the corresponding class
in ${\textsf{c}_{\rm top}}(P')$.

It is easy to check that the map
\begin{gather*}
{\textsf{c}_{\rm top}}\colon \  {\textsf{STB}}_\psi(W) \to \textrm{topological $G$-spaces over }W,\\
\hphantom{{\textsf{c}_{\rm top}}\colon} \
(P\xrightarrow{\varphi}P')\mapsto
\big({\textsf{c}_{\rm top}}(P)\xrightarrow{{\textsf{c}_{\rm top}}(\varphi)}{\textsf{c}_{\rm top}}(P')\big)
\end{gather*}
is a functor that commutes with restrictions to open subsets of $W$.
\end{proof}

\begin{proof}[Details of Step 3]
  We start by extending the symplectic reduction theorem of
  Marsden--Wein\-stein and Meyer~\cite{MW, Meyer} to manifolds with corners.

\begin{Theorem}\label{thm:red-corner}
  Suppose $(M,\sigma)$ is a symplectic manifold with corners with a proper
  Hamiltonian action of a Lie group $K$ and an associated equivariant
  moment map $\Phi\colon M\to{\mathfrak k}^*$.  Suppose further:
\begin{enumerate}[$1.$]\itemsep=0pt
\item For any point $x\in \Phi^{-1} (0)$ the stabilizer $K_x$ of $x$ is trivial;
\item there is an extension $\tilde{\Phi}$ of $\Phi$ to a manifold
  $\tilde{M}$ containing $M$ as a domain $($q.v.\
  Definition~{\rm \ref{def:domain})} with $\Phi^{-1} (0) = \tilde{\Phi}^{-1}
  (0)$.
\end{enumerate}
Then $\Phi^{-1} (0)$ is a manifold $($without corners$)$ and the quotient
\begin{gather*}
M/\!/_0 K:= \Phi^{-1} (0)/K
\end{gather*}
is naturally a symplectic manifold.
\end{Theorem}
\begin{Remark}
  The main issue in proving the theorem is in showing that $\Phi^{-1}
  (0)$ is actually a manifold and that it has the right dimension.  In
  other words the issue is transversality for manifolds with corners.
  To be more specif\/ic if $Q$ is a manifold with corners, $f\colon Q\to {\mathbb R}^k$
  is a~smooth function and $0$ is a regular value of $f$, then it is
  not true in general that $f^{-1} (0)$ is a~manifold, with or without
  corners. Take, for example,
\begin{gather*}
Q =\big\{(x,y,z)\in {\mathbb R}^3 \,|\, z\geq 0\big\}
\end{gather*}
and $f(x,y,z) = z-x^2 + y^2$.  Then $0$ is a regular value of $f$ but{\samepage
\begin{gather*}
f^{-1} (0) = \big\{(x,y,z)\in {\mathbb R}^3 \,|\,  z= x^2 -y^2,\,  z\geq 0\big\},
\end{gather*}
which is clearly not a manifold, with or without boundary.}

The standard approach to transversality for manifolds with
corners \cite{Nielsen} is to impose an additional requirement that the
kernel of the dif\/ferential of $f$ is transverse to the strata of $Q$.
However, in the situation we care about we have tangency instead.
Moreover, it is easy to write down an example of a smooth function
$h\colon {\mathbb R}^2\to {\mathbb R}$ so that the graph of $h$ is tangent to the $x$--$y$ plane
but the set
\begin{gather*}
\{(x,y,z)\,|\, z-h(x,y)=0,\, z\geq 0\}
\end{gather*}
is not a manifold.  This is why we make an awkward assumption on the
level set $\Phi^{-1} (0)$ in Theorem~\ref{thm:red-corner}.  On the other
hand, this assumption is easy to check in practice.
\end{Remark}

Before proving the theorem we f\/irst prove
\begin{Lemma}\label{lem2:28}
Let $f\colon Q\to {\mathbb R}^n$ be a smooth function on a manifold with corners~$Q$.
Suppose $\tilde{Q}$ is a manifold $($without corners$)$ containing~$Q$ as a domain,
and $\tilde{f}\colon \tilde{Q}\to {\mathbb R}^n$ is an extension of~$f$ with
\begin{gather*}
f^{-1} (0) = \tilde{f}^{-1} (0).
\end{gather*}
If $0$ is a regular value of $f$ $($that is, if for all $x\in f^{-1} (0)$
the map $d_x f\colon T_xQ\to {\mathbb R}^n$ is onto$)$, then~$f^{-1} (0)$ is naturally
a smooth manifold of dimension $\dim Q-n$
in the sense of Definition~{\rm \ref{submanifold}}.
\end{Lemma}

\begin{proof}
Since $0$ is a regular value of~$f$, and since $f^{-1} (0) = \tilde{f}^{-1} (0)$,
the value $0$ is also regular for~$\tilde{f}$.  Consequently,
$\tilde{f}^{-1} (0)$ is naturally a manifold of dimension $\dim Q - n$.
Since $\tilde{f}^{-1} (0) = f^{-1} (0)\subset Q$, we conclude
that~$f^{-1} (0)$ is naturally a manifold.
\end{proof}

\begin{Remark}
  Note that the assumptions of the lemma force $\ker d_x f$ to be
  tangent to the strata of $Q$: otherwise $f^{-1} (0) = \tilde{f}^{-1}
  (0)$ cannot hold.
\end{Remark}

\begin{proof}[Proof of Theorem~\ref{thm:red-corner}]
  Once we know that $\Phi^{-1} (0)$ is actually a manifold, the
  classical arguments of Marsden--Weinstein~\cite{MW} and of Meyer~\cite{Meyer} apply to show that $\sigma|_{\Phi^{-1} (0)}$ is basic and
  that its kernel is precisely the directions of the $G$ orbits.
  Consequently the restriction $\sigma|_{\Phi^{-1} (0)}$ descends to a
  closed nondegenerate 2-form $\sigma_0$ on the manifold $\Phi^{-1}
  (0)/K$.

  By Lemma~\ref{lem2:28} it is enough to show that 0 is a regular
  value of $\Phi$. This will follow from our assumption that the $K$
  action on $\Phi^{-1} (0)$ is free. Again, the argument is
  standard. Indeed, let $x \in \Phi^{-1} (0)$.  To show that the
  dif\/ferential $d_x \Phi\colon  T_xM \to {\mathfrak k}^*$ is surjective, we need to
  show that the annihilator of its image is zero.  Let $X \in {\mathfrak k}$ be
  in the annihilator of this image:
\begin{gather*}
\langle d_x \Phi (v), X\rangle = 0
      \qquad  \textrm{for all} \ \ v\in T_x M  .
\end{gather*}
By the def\/inition of the moment map, we may rewrite this as
\begin{gather*}
     \sigma_x(v,X_M(x)) = 0 \qquad \text{for all} \ \ v \in T_xM .
\end{gather*}
Since $\sigma_x$ is a nondegenerate form on $T_xM$,
we conclude that $X_M(x) = 0$.
Because the stabilizer of~$x$ is trivial, this implies
that $X$ is the zero vector in~${\mathfrak k}$.
\end{proof}

\begin{Remark}\label{rmrk:2.28}
  If additionally there is a Hamiltonian action of a Lie group $G$ on
  $(M,\sigma)$ with a moment map $\mu\colon M\to {\mathfrak g}^*$ so that the actions
  of~$G$ and~$K$ commute and the moment map~$\mu$ is~$K$ invariant
  and~$\Phi$ is $G$ invariant then
\begin{itemize}\itemsep=0pt
\item the induced action of $G$ on $(M/\!/_0 K, \sigma_0)$ is
  Hamiltonian and
\item $\mu|_{\Phi^{-1} (0)}$ descends to a moment map
  $\tilde{\mu}\colon M/\!/_0 K \to {\mathfrak g}^*$ for the induced action of $G$.
\end{itemize}
\end{Remark}

\begin{Lemma}
  Let $\psi\colon W\to {\mathfrak g}^*$ be a u.l.e.\ and $(\pi\colon P\to W,\omega)$ a
  symplectic toric $G$-bundle.  Then for every point $w\in W$ there is a~neighborhood $U_w$ so that the symplectic quotient~$(P|_{U_w} \times
  V_w)/\!/_0 K_w$ is a symplectic toric
  $G$-manifold.
\end{Lemma}

\begin{proof}
By Step~1  we have a
  neighborhood $U$ of $w\in W$ with faces, a subtorus $K=K_w$ of $G$, a basis
  $\{v_1,\ldots v_k\}$ of the integral lattice of ${\mathbb Z}_K$, the dual
  basis $\{v_1^*,\ldots, v_k^*\}$ of the weight lattice and  a~symplectic representation $K\to \operatorname{Sp}(V_w, \omega_w)$ with the weights
  $\{v_1^*,\ldots, v_k^*\}$
  so that the map
\begin{gather*}
\psi|_{U}\colon \ U\to C_w:= \{ \eta \in{\mathfrak g}^* \,|\, \langle \eta - \psi(w),
v_i\rangle \geq 0, \, 1\leq i\leq k\}
\end{gather*}
is an open embedding of manifolds with
corners.
It will be convenient to take $V_w= {\mathbb C}^k$, $\omega_w=
\frac{\sqrt{-1}}{2\pi}\sum dz_j \wedge d\bar{z}_j$ with the  action
of $K$ on ${\mathbb C}^k$ given by
\begin{gather*}
\exp (X)\cdot z:= \big(e^{2\pi \sqrt{-1}\langle v_1^*, X\rangle} z_1, \ldots,
e^{2\pi \sqrt{-1}\langle v_k^*, X\rangle} z_k\big).
\end{gather*}
We may do so by Lemma~\ref{lemma:step1.2}.
Then the moment map $\mu_w\colon {\mathbb C}^k\to {\mathfrak k}^*$ is given by the formula
\begin{gather*}
\mu_w (z) = -\sum |z_j|^2 v_j^*.
\end{gather*}
Let $\iota\colon {\mathfrak k}\hookrightarrow {\mathfrak g}$ denote the canonical inclusion and
$\iota^*\colon  {\mathfrak g}^*\to {\mathfrak k}^*$ the dual map.  Note that the kernel of $\iota^*$
is the annihilator ${\mathfrak k}^\circ $ of ${\mathfrak k}$ in ${\mathfrak g}^*$. Set
\begin{gather*}
\xi_0: = \iota^* (\psi (w)).
\end{gather*}
The cone
\begin{gather*}
C_w':= \big\{ \xi \in{\mathfrak k}^* \,|\, \langle \xi - \xi_0,
v_i\rangle \geq 0, \, 1\leq i\leq k\big\}
\end{gather*}
contains no nontrivial af\/f\/ine subspaces.
If we identify ${\mathfrak k}^*$ with a subspace of ${\mathfrak g}^*$,
the cone $C_w$ becomes the product ${\mathfrak k}^\circ \times C_w'$.
Since~$\psi|_{U}$ is an open embedding, we may assume (by shrinking~$U$
further if necessary) that $U$ is a product:
\begin{gather*}
U= {\mathcal O}\times U',
\end{gather*}
where ${\mathcal O}$ is a neighborhood of $w$ in the stratum containing it,
 $U' = {\mathcal V} \cap C'_w$
with ${\mathcal V}$ a neighborhood in ${\mathfrak k}^*$ of the apex of the cone $C'_w$,
and such that ${\mathcal O}$ and $U'$ are contractible.
{\it We take $U_w$ to be this~$U$.}
 Then, since $U$ is contractible
the restriction of any principal bundle $P\to W$ to $U$ is trivial.

Let $\nu\colon P|_U\to {\mathcal V}\cap C'_w \subset {\mathfrak k}^*$ be the composite
\begin{gather*}
\nu:= \iota^*\circ \psi \circ \pi.
\end{gather*}
We observe that~(1) $\nu$ is a moment map for the action of $K$ on
$(P|_U, \omega)$ and (2) $\nu\colon P|_U\to {\mathcal V}\cap C'_w$ is a trivial
f\/iber bundle (the typical f\/iber is ${\mathcal O}\times G$). Observation~(2)
implies that $\nu$ can be extended to a trivial ${\mathcal O}\times G$ f\/iber
bundle $\tilde{\nu}\colon \tilde{P} \to {\mathcal V}$ so that $P|_U$ embeds into
$\tilde{P}$ as a domain.  Observation (1) tells us that the diagonal
action of $K$ on $(P|_U\times {\mathbb C}^k, \omega\oplus \omega_w)$ is
Hamiltonian with a corresponding moment map
$\Phi \colon  P|_U \times {\mathbb C}^k \to {\mathfrak k}^*$ given by
\begin{gather*}
\Phi(p,z) = \nu(p) -\xi_0+ \mu_w (z) = \nu(p) -\xi_0 -\sum |z_j|^2 v_j^*.
\end{gather*}
Clearly
\begin{gather*}
\tilde{\Phi}(p,z) := \tilde{\nu}(p) -\xi_0 + \mu_w (z)
\end{gather*}
is an extension of $\Phi$. Since
\begin{gather*}
\tilde{\Phi}^{-1} (0) = \{(p,z)\,|\, \tilde{\nu}(p) = \xi_0-\mu_w(z)\}
\end{gather*}
and since
\begin{gather*}
\xi_0-\mu_w(z)\in C'_w \qquad \textrm{for all} \ \ z\in {\mathbb C}^k,
\end{gather*}
we have
\begin{gather*}
\Phi^{-1} (0) = \tilde{\Phi}^{-1} (0).
\end{gather*}
Since the action of $G$ on $P$ is free, so is the action of $K$ on
$P|_U\times {\mathbb C}^k$.  Therefore we can apply
Theorem~\ref{thm:red-corner} and conclude that
\begin{gather*}
\operatorname{cut}(P|_U):= (P|_U\times {\mathbb C}^k)/\!/_0 K
\end{gather*}
is a symplectic manifold (without corners).

The action of $G$ on $P$ extends trivially to a Hamiltonian action of~$G$ on $P|_U\times {\mathbb C}^k$.  This action of $G$ is Hamiltonian, commutes
with the action of $K$ and satisf\/ies the rest of the conditions of
Remark~\ref{rmrk:2.28}.  Consequently $\operatorname{cut}(P|_U)$ is a Hamiltonian
$G$-space.  Note that
\begin{gather*}
\dim \operatorname{cut}(P|_U) = \dim \big(P|_U \times {\mathbb C}^k\big) - 2k = \dim P = 2\dim G.
\end{gather*}
Thus to show that $\operatorname{cut}(P|_U)$ is toric, it is enough to show
that the action of $G$ is free at some point of $\operatorname{cut}(P|_U) $.  Now
take any point $\xi \in {\mathcal V}$ that also lies in the interior of the
cone~$C'_w$.  Pick any point $p\in P|_U$ with $\nu(p) = \xi$ and $z\in
{\mathbb C}^k$ with $\mu_w(z) = -\xi +\xi_0$.  Then
\begin{gather*}
\Phi(p,z) = \nu(p)+\mu_w(z) = \xi -\xi_0 +(-\xi +\xi_0) =0.
\end{gather*}
  On the other hand, since
$\xi$ is in the interior of the cone, the stabilizer of $z$ is
trivial.  Hence the stabilizer of $(p,z) \in P|_U\times {\mathbb C}^k$ for the
action of $G\times K$ is trivial as well. Consequently the stabilizer
of the image of $(p,z)$ in $\operatorname{cut}(P|_U)$ for the action of $G$ is
trivial.

We leave it to the reader to check that $\operatorname{cut}(P|_U)$ is a toric
manifold over $\psi|_U\colon U\to {\mathfrak g}^*$.
\end{proof}

\begin{Remark}
  If $(\pi_i\colon P_i\to W, \omega_i)$, for $i=1,2$, are two symplectic toric
  $G$-bundles over a u.l.e.\ $\psi\colon W\to {\mathfrak g}^*$ and
  $\varphi\colon P_1\to P_2$ is a morphism in ${\textsf{STB}}_\psi(W)$, i.e., a~$G$-equivariant symplectomorphism with $\pi_2 \circ \varphi =
  \pi_1$, then for any $w\in W$
\begin{gather*}
\varphi \times \text{id}\colon \  P_1|_{U_w}\times {\mathbb C}^k \to P_2|_{U_w}\times {\mathbb C}^k
\end{gather*}
is a $G\times K$-equivariant symplectomorphism with $\Phi_2\circ
(\varphi \times \text{id}) = \Phi_1$.
Hence $\varphi \times \text{id}$ maps $\Phi_1^{-1} (0)$
onto $\Phi_2^{-1} (0)$ and descends to an isomorphism of toric
manifolds
\begin{gather*}
\operatorname{cut}(\varphi)\colon \ \operatorname{cut}(P_1|_{U_w}) \to \operatorname{cut}(P_2|_{U_w}).
\end{gather*}
It is not hard to check that
\begin{gather*}
\operatorname{cut}\colon \  {\textsf{STB}}_\psi (U_w)\to {\textsf{STM}}_\psi (U_w)
\end{gather*}
is a functor for every $w\in W$.  (Strictly speaking we have a family of
functors parameterized by the points $w$ of $W$; we suppress this
dependence in our notation.)
\end{Remark}

We now proceed to construct the natural $G$-equivariant homeomorphisms
\begin{gather*}
\alpha_w^P\colon \  {\textsf{c}_{\rm top}}(P|_{U_w}) \to \operatorname{cut}(P|_{U_w}).
\end{gather*}
The construction depends on the fact that $({\mathbb C}^k, \omega_w, \mu_w)$ is
a symplectic toric $K$-manifold over the cone $\mu_w ({\mathbb C}^k) = \{\eta
\in {\mathfrak k}^*\,|\, \langle \eta, v_i\rangle \leq 0 \text{ for } 1\leq i \leq k\}$.
Moreover,
\begin{enumerate}[1)]\itemsep=0pt
\item the map $\mu_w\colon {\mathbb C}^k \to \mu_w({\mathbb C}^k) $ has a continuous
  (Lagrangian) section
\begin{gather*}
s\colon \ \mu_w({\mathbb C}^k) \to {\mathbb C}^k, \qquad s(\eta) =
\big(\sqrt{\langle -\eta, v_1\rangle},\ldots, \sqrt{\langle -\eta, v_k\rangle}\big)
\end{gather*}
which is smooth over the interior of the cone $\mu_w ({\mathbb C}^k)$;
\item the stabilizer $K_z$ of $z\in {\mathbb C}^k$ depends only on the face of
  the cone $\mu_w ({\mathbb C}^k)$ containing $\mu_w(z)$ in its interior:
\begin{gather*}
K_z = \exp \big(\operatorname{span}_{\mathbb R} \{v_i \in \{v_1,\ldots, v_k\} \,|\, \langle
  \mu_w(z) , v_i\rangle = 0\}\big);
\end{gather*}
cf.\ Remark~\ref{rmrk:2.21}.
\end{enumerate}
We continue with the notation above: $\xi_0 = \iota^*(\psi(w)) \in
{\mathfrak k}^*$ is a point and $\nu= \iota^*\circ \mu \colon P|_U\to {\mathfrak k}^*$ the
$K$-moment map.  Then for any point $p\in P|_U$
\begin{gather*}
\xi_0 -\nu(p) \in \mu_w\big({\mathbb C}^k\big)
\end{gather*}
and
\begin{gather*}
s(\xi_0 -\nu(p) )  =  \big(\sqrt{\langle \nu(p) -\xi_0, v_1\rangle},\ldots,
\sqrt{\langle \nu(p) -\xi_0, v_k\rangle}\big)\\
\hphantom{s(\xi_0 -\nu(p) )}{}
 =  \big(\sqrt{\langle \mu(p) -\psi (w),
  v_1\rangle},\ldots, \sqrt{\langle \mu(p) -\psi(w), v_k\rangle}\big),
\end{gather*}
where $\mu = \psi \circ \pi\colon P\to {\mathfrak g}^*$ is the moment map for the
action of $G$ on $P$.  This gives us a~continuous proper map
\begin{gather*}
\phi\colon \  P|_U\to \Phi^{-1} (0) \subset P|_U \times {\mathbb C}^k, \qquad
     \phi(p) = (p, s (\xi_0 - \nu(p))).
\end{gather*}
The image of $\phi$
intersects every $K$ orbit in $\Phi^{-1} (0)$.  Hence the composite
\begin{gather*}
f = \tau\circ \phi\colon \ P|_U\to \Phi^{-1} (0)/K,
\end{gather*}
where $\tau\colon  \Phi^{-1} (0)\to \Phi^{-1} (0)/K$ is the orbit map, is
surjective.  Next we argue that the f\/ibers of $f$ are precisely
the equivalence classes of the relation $\sim$ def\/ined in Step~2.
Two points $p_1, p_2\in P|_U$ are equivalent
with respect to $\sim$ if and only if $\pi (p_1) = \pi(p_2)$ and there
is an $a\in K_{\pi (p_1)}$ with $a\cdot p_2 = p_1$.  On the other hand
$f(p_1) = f(p_2)$ if and only if there is an~$a\in K$ with
\begin{gather*}
(p_1,  s (\xi_0 - \nu(p_1))) = (a\cdot p_2, a\cdot s(\xi_0 - \nu(p_2))).
\end{gather*}
For any point $\eta\in \mu_w({\mathbb C}^k)$
\begin{gather*}
a\cdot s(\eta) = s(\eta) \quad \Leftrightarrow \quad a
  \textrm{ lies in the stabilizer }K_{s(\eta)}\textrm{ of } s(\eta).
\end{gather*}
For $\eta = \xi_0 -\nu(p_2) =\xi_0 -\nu(p_1)
 =  \iota^*(\psi(w)-\psi (\pi(p_1)))$,
\begin{gather*}
K_{s(\eta)}  =  \exp \big(\operatorname{span}_{\mathbb R} \{v_i \in \{v_1,\ldots, v_k\} \,|\, \langle
  \xi_0 - \nu(p_1), v_i\rangle = 0\}\big)\\
 \hphantom{K_{s(\eta)}}{}
 =  \exp \big(\operatorname{span}_{\mathbb R} \{v_i \in \{v_1,\ldots, v_k\} \,|\, \langle
  \iota^*(\psi(w)- \psi (\pi(p_1))) , v_i\rangle = 0\}\big)\\
\hphantom{K_{s(\eta)}}{} =  \exp \big(\operatorname{span}_{\mathbb R} \{v_i \in \{v_1,\ldots, v_k\} \,|\, \langle
  \psi(w)- \psi (\pi(p_1)), v_i\rangle = 0\}\big) \\
\hphantom{K_{s(\eta)}}{} =  K_{\pi(p_1)}.
\end{gather*}
We conclude that the f\/ibers of $f$ are precisely
the equivalence classes of the relation $\sim$.
Therefo\-re~$f$ descends to a continuous bijection
\begin{gather*}
\alpha^P_w\colon  \  {\textsf{c}_{\rm top}} (P|_U)=(P|_U)/_{\sim}   \to
   \Phi^{-1} (0)/K= \operatorname{cut}(P|_U), \qquad \alpha^P_w([p])=[p,s(\xi_0 - \nu(p))] .
\end{gather*}
The properness of $f$ implies that $\alpha^P_w$ is a homeomorphism.
This follows from Lemma~\ref{compactly generated} below.

\begin{Lemma} \label{compactly generated}
Let $f \colon A \to B$ be a continuous and proper bijection
between topological spaces.  Suppose that $B$ is Hausdorff
and compactly generated.  That is, $B$ is Hausdorff and a subset~$E$ of~$B$ is closed if and only if for every compact~$K$
the intersection~$E \cap K$ is compact.
Then~$f$ is a~homeomorphism.
\end{Lemma}

\begin{proof}
Omitted.
\end{proof}

The commutativity of \eqref{eq:alpha-cut} is easy:
Since $\nu_1(p) = \nu_2(\varphi(p))$,
\begin{gather*}
  \operatorname{cut}(\varphi) \big(\alpha^{P_1}_w [p]\big) =  \operatorname{cut}(\varphi) ([p,s(\xi_0 - \nu_1(p))])
    =  [\varphi(p),s(\xi_0 - \nu_2(\varphi(p)))]\\
\hphantom{\operatorname{cut}(\varphi) \big(\alpha^{P_1}_w [p]\big)}{} = \alpha^{P_2}_w ([\varphi(p)])
  =    \alpha^{P_2}_w ({\textsf{c}_{\rm top}}(\varphi) ([p])).
\end{gather*}

To f\/inish the construction of the functor ${\textsf{c}}$ it remains to show
that $\upsilon:= (\alpha^P_{w_2 })\circ (\alpha^P_{w_1})^{-1}$ is a~map
of symplectic toric $G$-manifolds.
Since $\alpha^P_{w_1}$ and $\alpha^P_{w_2}$
are $G$-equivariant homeomorphisms, so is $\upsilon$.
It is enough to produce a smooth symplectic map $\vartheta$ satisfying
\begin{gather*}
\vartheta \circ \alpha^P_{w_1} = \alpha^P_{w_2}.
\end{gather*}
Indeed, this implies that $\vartheta = \upsilon$,
and hence that $\upsilon$ is a smooth symplectic map;
reversing the roles of~$w_1$ and~$w_2$, we conclude that the inverse
of~$\upsilon$ is also a smooth symplectic map,
and so~$\upsilon$ is a $G$-equivariant dif\/feomorphism.
Since the intersection $U_{w_1} \cap U_{w_2}$ can be covered by sets
of the form~$U_{w_3}$, it suf\/f\/ices to consider the case when~$U_{w_1}$
is contained in~$U_{w_2}$.

Consider f\/irst the special case when $K_{w_1} = K_{w_2} =K$.  Then the
collections of the correspon\-ding weights $\{
v_j^{(w_1)}\}_{j=1}^k$, $\{ v_j^{(w_2)}\}_{j=1}^k$ are the same set.
Hence by Lemma~\ref{lemma:step1.2} there exists a~symplectic linear
isomorphism $\vartheta\colon {\mathbb C}^k\to {\mathbb C}^k$ which permutes coordinates and
intertwines the two representations and the corresponding moment maps.
Consequently ${\rm id}\times \vartheta\colon P|_{U_{w_1}}\times {\mathbb C}^k \to
  P|_{U_{w_1}} \times {\mathbb C}^k$ induces a symplectic isomorphism of
    symplectic quotients
\begin{gather*}
\overline{\vartheta}\colon \ \big(P|_{U_{w_1}}\times {\mathbb C}^k\big)/\!/_0 K\to
  \big(P|_{U_{w_1}}\times {\mathbb C}^k\big)/\!/_0 K,\qquad [p,z]\mapsto [p,
    \vartheta(z)].
\end{gather*}
It is easy to check that $\overline{\vartheta}\circ \alpha^P_{w_1} =
\alpha^P_{w_2}$, hence $(\alpha^P_{w_2 })\circ (\alpha^P_{w_1})^{-1}$
is a symplectomorphism in {\em this} case.

More generally we have a strict inclusion $\{
v_j^{(w_1)}\}_{j=1}^{k_1}\subset \{ v_j^{(w_2)}\}_{j=1}^{k_2}$.  By
the discussion of the special case above, it is not a loss of generality
to assume that $v_j^{(w_1)} = v_j^{(w_2)}$ for all $1\leq j\leq k_1$.
We may then reduce the clutter in the notation by dropping the
superscripts $^{(w_1)}$ and $^{(w_2)}$ and setting $K_i:=K_{w_i}$, $i=1,2$.

By construction of the neighborhoods $U_{w_i}$ (q.v.\
Lemma~\ref{lemma:step1.1} and subsequent remarks) we have
\begin{itemize}\itemsep=0pt
\item $\langle\psi(w_1)-\psi(w_2), v_i\rangle =0$ for $i=1,\ldots, k_1$, and,
\item for all $w\in U_{w_1}$,
\begin{gather*}
\langle \psi(w) - \psi(w_2), v_i\rangle>0 \qquad \textrm{for} \ \ i=k_1+1,\ldots, k_2.
\end{gather*}
\end{itemize}
Consequently  for any point $p\in P|_{U_{w_1}}$ the functions
\begin{gather*}
p\mapsto \sqrt{\langle \mu(p) - \psi(w_2),v_i \rangle}
\end{gather*}
are smooth for $i=k_1+1,\ldots, k_2$.  Also, for  $p\in P|_{U_{w_1}}$
\begin{gather*}
\langle \mu(p) - \psi(w_2),v_i \rangle = \langle \mu(p) - \psi(w_1),v_i \rangle
\end{gather*}
for $i=1,\ldots, k_1$.  Now consider the map
\begin{gather*}
\vartheta\colon \ P|_{U_{w_1}}\times {\mathbb C}^{k_1} \to P|_{U_{w_2}} \times {\mathbb C}^{k_2}
\end{gather*}
given by
\begin{gather*}
\vartheta(p,z_1,\ldots, z_{k_1}) = \Big(p,z_1,\ldots z_{k_1},
\sqrt{\langle   \mu(p) - \psi(w_2),v_{k_1+1} \rangle},\ldots ,
\sqrt{\langle \mu(p) -  \psi(w_2),v_{k_2} \rangle}\Big).
\end{gather*}
The map $\vartheta $ is smooth and $K_1$-equivariant.  Since
$\vartheta^*(dz_j\wedge d\bar{z}_j) = 0$ for $j>k_1$, it is symplectic.

Next observe that
\begin{gather*}
\vartheta^{-1} \big(\Phi_2 ^{-1} (0)\big) = \Phi_1 ^{-1} (0),
\end{gather*}
where $\Phi_j\colon P|_{U_{w_1}}\times {\mathbb C}^{k_j}\to {\mathfrak k}_j^*$, $j=1,2$ are
the corresponding moment maps.  This is because
\begin{gather*}
(p,z)\in \Phi_j^{-1} (0)\Leftrightarrow
\langle \psi(\pi(p)) - \psi(w_j), v_i\rangle = |z_i|^2
\qquad\textrm{for all} \ \ i=1,\ldots, k_j .
\end{gather*}
Consequently $\vartheta$ descends to a well-def\/ined smooth symplectic map
\begin{gather*}
\bar{\vartheta}\colon \  \Phi_1^{-1} (0)/K_1 \to \Phi_2 ^{-1} (0)/K_2
\end{gather*}
given by
\begin{gather*}
\bar{\vartheta}([p, z_1,\ldots,z_{k_1}]) = \Big[p, z_1,\ldots,z_{k_1},
\sqrt{\langle   \mu(p) - \psi(w_2),v_{k_1+1} \rangle},\ldots ,
\sqrt{\langle \mu(p) -  \psi(w_2),v_{k_2} \rangle }\Big].
\end{gather*}
Evidently,
\begin{gather*}
\bar{\vartheta} \big(\alpha^P_{w_1}([p])\big) = \alpha^P_{w_2}([p]).
\end{gather*}
This f\/inishes Step~3 of the construction.  We have thus constructed
the desired functor~${\textsf{c}}$.
\end{proof}

We are now in position to prove part (1) of Theorem~\ref{thm:main}.

\begin{proof}[Proof of Theorem~\ref{thm:main}(1)]
  It is enough to show that the category ${\textsf{STB}}_\psi (W)$ is nonempty.
  For then for any object $P$ of ${\textsf{STB}}_\psi (W)$, the object ${\textsf{c}} (P)$
  is the desired symplectic toric manifold.

  Consider the cotangent bundle $T^*G$ with the action of $G$ given by
  the lift of multiplication on the left.  This action is Hamiltonian
  and the moment map $\mu\colon T^*G\to {\mathfrak g}^*$ makes $T^*G$ into a~symplectic principal $G$-bundle over ${\mathfrak g}^*$.  The pullback of this
  bundle by $\psi\colon W\to {\mathfrak g}^*$ is an object of ${\textsf{STB}}_\psi (W)$.

  Alternatively, for any principal $G$-bundle $\pi\colon P\to W$ and any
  choice of a connection 1-form $A\in \Omega^1(P,{\mathfrak g})^G$ the closed
  2-form $\sigma: = d\langle \psi \circ \pi ,A\rangle$ is symplectic
  and $\mu:= \psi \circ \pi$ is a corresponding moment map (see
  Lemma~\ref{lemma:nondeg} below). Then $(P\xrightarrow{\pi} W,
  \sigma)$ is an object of ${\textsf{STB}}_\psi (W)$.
\end{proof}

We end the section with a lemma that will be used to prove that
$\textsf{c}\colon {\textsf{STB}}_\psi(W)\to {\textsf{STM}}_\psi (W)$ is an equivalence of categories.

\begin{Lemma} \label{lem: c is a map of stacks} The functor
  ${\textsf{c}}\colon {\textsf{STB}}_\psi \to {\textsf{STM}}_\psi$ is a map of presheaves of groupoids.
  Moreover over the interior ${\mathaccent23{W}}$ the functor ${\textsf{c}}_{{\mathaccent23{W}}}\colon
  {\textsf{STB}}_\psi({\mathaccent23{W}}) \to {\textsf{STM}}_\psi({\mathaccent23{W}})$ is isomorphic to the
  identity functor.
\end{Lemma}

\begin{proof}
  Since ${\textsf{c}_{\rm top}}$ is a map of presheaves then so is $\textsf{c}$.  Over the
  interior ${\mathaccent23{W}}$ the functor ${\textsf{c}}_{\rm top}$ is isomorphic to the
  identity functor, since we divide out by the relation whose
  equivalence classes are singletons.

  Moreover, for any point $w\in {\mathaccent23{W}}$ the corresponding group $K_w$
  is trivial.  Hence
\begin{gather*}
\operatorname{cut}(P|_{U_w}) = (P|_{U_w}\times \{0\})/\!/_0 \{1\} \simeq P|_{U_w}
\end{gather*}
as symplectic toric manifolds.
\end{proof}

\section[Local trivializations of symplectic toric $G$-bundles]{Local trivializations of symplectic toric $\boldsymbol{G}$-bundles}
\label{sec:local uniqueness STB}

The purpose of this section is to prove the ``local uniqueness''
for symplectic toric $G$ bundles.  Here is the statement:

\begin{Lemma} \label{loc-iso-stb} Let $P_0 = ( \pi_0 \colon P_0 \to
  W,\omega_0)$ and $P_1 = (\pi_1 \colon P_1 \to W,\omega_1)$ be two
  symplectic toric $G$-bundles over a unimodular local embedding $($u.l.e.$)$
  $\psi \colon W \to {\mathfrak g}^*$.
  Then for any  open subset~$U$ of~$W$ with $H^2(U,{\mathbb Z}) =0$ the
  restrictions $P_1|_U$ and $P_1|_U$ are isomorphic in ${\textsf{STB}}_\psi (U)$.
Consequently, any two  symplectic principal $G$-bundles over the
same u.l.e.\ are locally isomorphic.
\end{Lemma}
Before proving Lemma~\ref{loc-iso-stb}, we need to establish
two facts about symplectic forms on principal  $G$-bundles over our
u.l.e.\ $\psi \colon W \to {\mathfrak g}^*$.  Recall that the notion of a moment map
does not {\em a priori} require a 2-form to be nondegenerate,
we only need~\eqref{eq:moment_map} to hold.

Note that it makes sense to pair a connection 1-form $A$ on $P$
with the moment map $\mu$.  This results in a real-valued $G$-invariant 1-form
  $\langle \mu, A\rangle\in \Omega^1 (P)^G$.

\begin{Lemma}\label{lemma:nondeg}
  Let $\psi\colon W\to {\mathfrak g}^*$ be a u.l.e.\ and $\pi\colon P \xrightarrow{\pi}W$
  a principal $G$-bundle.
\begin{itemize}\itemsep=0pt
\item Any closed $G$-invariant $2$-form on $P$ with moment map $\mu:=
  \psi \circ \pi$ is automatically symplectic.
\item A connection $1$-form $A\in \Omega^1(P,{\mathfrak g})^G$ defines a bijection
  from the space of closed forms on~$W$ and the space of invariant
  symplectic forms on~$P$ with moment map~$\mu$, by
\begin{gather*}
\beta \mapsto d \langle \mu, A\rangle + \pi^*\beta.
\end{gather*}
Consequently, any two closed $G$-invariant $2$-forms on $P$ with moment
map $\mu$ differ by a~basic closed $2$-form.
\end{itemize}
\end{Lemma}

\begin{proof}
  We argue f\/irst that $d \langle \mu, A\rangle$ is nondegenerate.

  For any vector $X\in {\mathfrak g}$ the Lie derivative $L_{X_P} \langle \mu,
  A\rangle$ with respect to the induced vector f\/ield~$X_P$ is zero.
  By Cartan's formula we then have
\begin{gather*}
0= d \iota(X_P)\langle  \mu, A\rangle
   +  \iota(X_P)d \langle  \mu, A\rangle =
   d \langle \mu, X\rangle + \iota(X_P)d \langle  \mu, A\rangle.
\end{gather*}
Therefore $\mu$ is a moment map for the action of~$G$
on $(P, d \langle \mu, A\rangle)$.  Also,
\begin{gather*}
d\langle \mu, A\rangle = \langle d\mu\wedge A\rangle + \langle \mu,
dA\rangle.
\end{gather*}
Moreover, for any point $p\in P$
we have an isomorphism
\begin{gather*}
T_pP ={\mathcal H}_p \oplus {\mathcal V}_p \xrightarrow{(d\mu_p\oplus A_p)} {\mathfrak g}^* \oplus {\mathfrak g},
\end{gather*}
where ${\mathcal H}_p$ and ${\mathcal V}_p$ are the horizontal and vertical subspace
of $T_pP$ respectively.
The map $d\mu_p|_{{\mathcal H}_p} \oplus
A_p|_{{\mathcal V}_p}$ is an isomorphism because $d\pi_P\colon {\mathcal H}_p \to
T_{\pi(p)}W$ is an isomorphism, $d\psi_{\pi(p)}\colon  T_{\pi(p)}W \to
{\mathfrak g}^*$ is an isomorphism since $\psi_p$ is an embedding,
and $A_p\colon {\mathcal V}_p\to {\mathfrak g}$ is an isomorphism too.  The isomorphism
$d\mu_p|_{{\mathcal H}_p} \oplus A_p|_{{\mathcal V}_p}$ identif\/ies $ \langle
d\mu\wedge A\rangle$ with the canonical pairing ${\mathfrak g}^*\times {\mathfrak g}\to
{\mathbb R}$. On the other hand $\langle \mu, dA\rangle$ is basic.  It follows
that $d\langle \mu, A\rangle $ is nondegenerate.

Similarly, for any closed form $\beta$ on $W$ the form
\begin{gather*}
\omega = \omega_{A,\beta}:= d\langle \mu, A\rangle +\pi^* \beta
\end{gather*}
is a closed, nondegenerate $G$-invariant form on $P$ and $\mu$ is a~moment map for the action of~$G$ on~$(P,\omega)$.

Finally, if $\sigma\in \Omega^2(P)^G$ is a closed $G$-invariant 2-form
with moment map $\mu$ and $A\in \Omega^1(P,{\mathfrak g})^G$ is a connection
1-form, then for any $X\in {\mathfrak g}$
\begin{gather*}
\iota(X_P)(\sigma - d\langle \mu, A\rangle) = -d \langle \mu, X\rangle
 + d \langle \mu, X\rangle = 0.
\end{gather*}
Hence $\sigma - d\langle \mu, A\rangle$ is basic and there is a 2-form
$\beta$ on $W$ with
\begin{gather*}
\sigma = d\langle \mu, A\rangle  + \pi^*\beta;
\end{gather*}
$\beta$ is necessarily closed.  Note that this also proves that $\sigma$ is nondegenerate.
\end{proof}

\begin{Lemma} \label{unique omega}
Let $\psi\colon W\to {\mathfrak g}^*$ be a u.l.e.\
  and $\pi \colon  P\xrightarrow{}W$ a principal $G$-bundle as above.
  Suppose that $\omega\in \Omega^2(P)^G$ is a closed $G$-invariant form with
  moment map $\mu$ and $\gamma\in \Omega^1 (W)$ a~$1$-form.  Then
  $(P,\omega)$ and $(P,\omega + \pi^* d\gamma)$ are isomorphic in
  ${\textsf{STB}}_\psi (W)$.  That is, there exists a gauge transformation
  $f\colon P\to P$ with $f^*(\omega + \pi^* d\gamma) = \omega$.
\end{Lemma}

\begin{proof}
We apply Moser's deformation method~\cite{moser}.  By Lemma~\ref{lemma:nondeg}
the forms
\begin{gather*}
\omega_t = \omega + t \pi^* d\gamma, \qquad t\in [0,1],
\end{gather*}
are symplectic. They have  the same moment map~$\mu$.
Let $X_t$ be the time-dependent vector f\/ield on $P$ that satisf\/ies
\begin{gather*}
\iota(X_t)\omega_t = - \pi^* \gamma.
\end{gather*}
Note that $X_t$ is $G$-invariant.
For every $\xi \in {\mathfrak g}$, we have
\begin{gather*}
\iota(X_t) d\langle \mu, \xi\rangle = -\omega_t(\xi_P,X_t) =
 - \iota(\xi_P) \pi^* \gamma =0.
\end{gather*}
Hence $d \mu ( X_t) = 0$. Since  $\mu = \psi \circ \pi$
and $\psi$ is a local embedding we have
\begin{gather*}
d \pi  (X_t) = 0.
\end{gather*}
That is, $X_t$ is tangent to the f\/ibers of $P\to W$, which
are tori.
Consequently we can integrate the vector f\/ield $X_t$ to obtain a~$G$-equivariant isotopy $\phi_t \colon P\to P$ which exists for all
$t\in [0,1]$
and projects to the identity map on the base $W$.  Then
\begin{gather*}
\frac{d}{dt} (\phi_t^* \omega_t)
 = \phi_t^* \left( L_{X_t} \omega_t + \frac{d}{dt} \omega_t\right)
 = d \phi_t^* (\iota(X_t) \omega_t + \pi^* \gamma) = 0.
\end{gather*}
Consequently $f := \phi_1\colon (P,\omega)\to (P,\omega + \pi^* d\gamma)$ is
an isomorphism of symplectic toric $G$-bundles over $\psi\colon W\to {\mathfrak g}^*$,
as desired.
\end{proof}

\begin{proof}[Proof of Lemma~\ref{loc-iso-stb}]\sloppy
Recall that for a torus $G= {\mathfrak g}/{\mathbb Z}_G$ the principal $G$-bundles
over a~mani\-fold with corners $N$  are classif\/ied by $H^2(N,{\mathbb Z}_G)$.
Since $H^2(U, {\mathbb Z}) = 0$ by assumption, $H^2(U,{\mathbb Z}_G) = 0$ as well.
  Consequently there exists a $G$-equivariant dif\/feomorphism
\begin{gather*}
h\colon \  P_0|_U \to P_1|_U
\end{gather*}
inducing the identity map on $U$.  By Lemma~\ref{lemma:nondeg}
\begin{gather*}
h^*\omega_1 = \omega_0 + \pi^*\beta
\end{gather*}
for some closed 2-form $\beta$ on $U$.  Since $H^2 (U, {\mathbb R}) =0$, there is a 1-form $\gamma$ on $U$ with $\beta = d \gamma$.  By Lemma~\ref{unique omega} there is a gauge transformation $f\colon  P_0|_U \to P_0|_U$ with
\begin{gather*}
f^* (\omega_0 + \pi^*\beta) = \omega_0.
\end{gather*}
Therefore
\begin{gather*}
(h\circ f)^* \omega_1 = f^* (h^*\omega_1) = f^* (\omega_0 +
\pi^*\beta) = \omega_0.\tag*{\qed}
\end{gather*}
\renewcommand{\qed}{}
\end{proof}

\begin{Remark}
  In the language of stacks Lemma~\ref{loc-iso-stb} asserts that the
  stack ${\textsf{STB}}_\psi$ is a gerbe: any two objects are locally
  isomorphic.
\end{Remark}

\section{Equivalence of categories of symplectic toric bundles\\
and symplectic toric manifolds}
\label{sec:equiv of cat}

In this section we show that the functor ${\textsf{c}}$
is an equivalence of categories.
This reduces the classif\/ication of symplectic toric $G$-manifolds
to that of symplectic toric $G$-bundles.  More specif\/ically we prove

\begin{Theorem}\label{thm:iso}
Let $\psi\colon W\to{\mathfrak g}^*$ be a unimodular local embedding $($u.l.e$)$.  The functor
\begin{gather*}
{\textsf{c}}\colon \  {\textsf{STB}}_\psi(W)\to {\textsf{STM}}_\psi (W)
\end{gather*}
constructed in Section~{\rm \ref{sec:collapse}} is an equivalence of
categories.
\end{Theorem}
\begin{Remark}
We actually show that for any open set $U\subset W$ the functor
\begin{gather*}
{\textsf{c}}_U\colon \  {\textsf{STB}}_\psi(U)\to {\textsf{STM}}_\psi (U)
\end{gather*}
is an equivalence of categories.  In other words $c$ is an isomorphism
of presheaves of groupoids.
\end{Remark}

  The proof of Theorem~\ref{thm:iso} proceeds in a series of lemmas.
Our f\/irst goal is to
  prove that the functor ${\textsf{c}}$ is full and faithful.

\begin{Lemma}[${\textsf{c}}$ is faithful] \label{lem4:1}
  For any open subset $U$ of $W$ and for any two objects $P_1, P_2 \in
  {\textsf{STB}}_\psi (U)$ the map
\begin{gather*}
{\textsf{c}}= {\textsf{c}}_U \colon  \ \operatorname{Hom}(P_1, P_2) \to \operatorname{Hom}({\textsf{c}}(P_1) {\textsf{c}}(P_2)), \qquad \phi\mapsto {\textsf{c}}(\phi)
\end{gather*}
is injective.
\end{Lemma}

\begin{proof}
The idea is easy: if two isomorphisms of symplectic toric $G$ bundles over~$W$
map to the same isomorphism of symplectic toric $G$ manifolds over~$W$,
then they must coincide over the interior of~$W$. By continuity,
they must coincide over all of~$W$.

  (a)~Recall that the functor
\begin{gather*}
c_{U\cap {\mathaccent23{W}}}\colon \  {\textsf{STB}}_\psi (U\cap {\mathaccent23{W}}) \to {\textsf{STM}}_\psi (U\cap {\mathaccent23{W}})
\end{gather*}
is isomorphic to the identity functor
(q.v.\ Lemma~\ref{lem: c is a map of stacks}):
we have isomorphisms
$\{ \delta_Q \colon  Q \to {\textsf{c}}(Q) \}_{Q\in {\textsf{STB}}_\psi (U\cap {\mathaccent23{W}})}$
such that for any morphism
$\phi\colon Q_1\to Q_2$ in ${\textsf{STB}} _\psi (U\cap {\mathaccent23{W}})$ the diagram
\begin{gather*}
\xy
(-20,10)*+{Q_1 }="1";
(20,10)*+ {{\textsf{c}}(Q_1)}="2";
(-20,-7)*+{Q_2}="3";
(20,-7)*+{{\textsf{c}}(Q_2)}="4";
 {\ar@{->}^{\delta_{Q_1}} "1";"2"};
{\ar@{->}_{{\textsf{c}}(\phi)}  "2";"4"};
 {\ar@{->}_{\phi} "1";"3"};
{\ar@{->}^{\delta_{Q_2}} "3";"4"};
\endxy
\end{gather*}
commutes.   Hence
\begin{gather*}
{\textsf{c}}\colon \  \operatorname{Hom}(Q_1, Q_2) \to  \operatorname{Hom}({\textsf{c}}(Q_1), {\textsf{c}}(Q_2))
\end{gather*}
is invertible with the inverse ${\textsf{c}}^{-1} $ given by
\begin{gather}\label{eq:cinv}
{\textsf{c}}^{-1} (\varphi) =\delta_{Q_2}^{-1} \circ \varphi \circ \delta_{Q_1}.
\end{gather}

 (b)~If $\phi_1,\phi_2\in \operatorname{Hom}(P_1, P_2)$ are two
morphisms with ${\textsf{c}}(\phi_1) = {\textsf{c}}(\phi_2)$ then
\begin{gather*}
{\textsf{c}}\big(\phi_1|_{P_1|_{U\cap {\mathaccent23{W}}}}\big) ={\textsf{c}}(\phi_1)|_{P_1|_{U\cap {\mathaccent23{W}}}} =
{\textsf{c}}(\phi_2)|_{P_1|_{U\cap {\mathaccent23{W}}}}={\textsf{c}}\big(\phi_2|_{P_1|_{U\cap {\mathaccent23{W}}}}\big).
\end{gather*}
By (b) above
\begin{gather*}
\phi_1|_{P_1|_{U\cap {\mathaccent23{W}}}}= \phi_2|_{P_1|_{U\cap {\mathaccent23{W}}}}.
\end{gather*}
Since the
restriction $P_1|_{U\cap {\mathaccent23{W}}}$ is dense in $P_1|_U$,
\begin{gather*}
\phi_1 = \phi_2.\tag*{\qed}
\end{gather*}
\renewcommand{\qed}{}
\end{proof}

Next, we need to check that, for every two objects $P_1$ and $P_2$,
every morphism $\varphi \colon {\textsf{c}}(P_1) \to {\textsf{c}}(P_2)$
comes from a morphism $P_1 \to P_2$.
Again, the idea is easy: $\varphi$ gives a morphism ${\mathaccent23{\phi}}$
between open dense subsets of~$P_1$ and~$P_2$,
(namely, the preimages of the interior of~$W$),
and we need to check that ${\mathaccent23{\phi}}$ extends smoothly to the boundary.
It is enough to check that ${\mathaccent23{\phi}}$ extends locally;
local extensions will coincide on the overlaps of their domains.
Locally, $P_1$ and $P_2$ are isomorphic, so it remains to consider the
case that $P_1 = P_2$.  For this case we will use
the following theorem of Haef\/liger, Salem and Schwartz:

\begin{Theorem}[\protect{\cite[Theorem~3.1]{HS}}] \label{thm:HSS}
Let
  $M$ be a manifold with an action of a torus $G$ and $h\colon M\to M$ a
  $G$-equivariant diffeomorphism with $h(x)\in G\cdot x$ for all
  points $x\in M$. Let $\pi\colon  M \to M/G$ be the orbit map.
  Then there exists a map $f\colon M/G\to G$ such that
\begin{gather*}
h(x) = f(\pi (x))\cdot x
\end{gather*}
for all $x\in M$ and such that $f \circ \pi$ is smooth.
\end{Theorem}

We continue with the proof that the functor ${\textsf{c}}$ is full.

\begin{Lemma} \label{lem4:2}
For any open subset $U$ of $W$ and for any $P\in {\textsf{STB}}_\psi(U)$ the map
\begin{gather*}
{\textsf{c}}\colon \ \operatorname{Hom}(P,P)\to \operatorname{Hom} ({\textsf{c}}(P), {\textsf{c}}(P))
\end{gather*}
is onto.
\end{Lemma}

\begin{proof}
  By Theorem~\ref{thm:HSS}, given $\varphi\in \operatorname{Hom}({\textsf{c}}(P),{\textsf{c}}(P))$ there
  is a smooth function $f\colon U\to G$ so that
\begin{gather*}
\varphi(x) = f(\pi(x))\cdot x,
\end{gather*}
where $\pi\colon {\textsf{c}}(P)\to U$ is the quotient map.  By Step~(a) of the proof
of Lemma~\ref{lem4:1} and~\eqref{eq:cinv},
\begin{gather*}
\varphi|_{P|_{U\cap {\mathaccent23{W}}}} = {\textsf{c}} ({\mathaccent23{\phi}}),
\end{gather*}
where ${\mathaccent23{\phi}}$ is given by
\begin{gather*}
{\mathaccent23{\phi}}= (\delta_P)^{-1} \circ \varphi|_{P|_{U\cap {\mathaccent23{W}}}} \circ \delta_P.
\end{gather*}
Hence for $p\in P|_{U\cap {\mathaccent23{W}}}$,
\begin{gather*}
{\mathaccent23{\phi}}(p) = (\delta_P)^{-1} \left(f (\pi(\delta_P (p)))\cdot
  \delta_P(p)\right)= (\delta_P)^{-1} \left(f(\pi(p))\cdot
  \delta_P(p)\right)\\
\hphantom{{\mathaccent23{\phi}}(p)}{}
= f(\pi(p))\cdot \delta_P^{-1} (\delta_P(p)) =
f(\pi(p))\cdot p.
\end{gather*}
Def\/ine the map $\phi\colon P\to P$ by
\begin{gather*}
\phi(p) = f(\pi(p))\cdot p\quad \textrm{ for all }p\in P.
\end{gather*}
This map is $G$-equivariant and commutes with the orbit map $\pi\colon P\to
U$.  Since $f \circ \pi$ is smooth, the map $\phi$ is a
dif\/feomorphism.  Moreover since the restriction of $\phi$ to $P|_{U\cap {\mathaccent23{W}}}$ is
${\mathaccent23{\phi}}$, the map $\phi$ is symplectic on $P|_{U\cap {\mathaccent23{W}}}$.  Since
$P|_{U\cap {\mathaccent23{W}}}$ is dense in $P$, we conclude that $\phi$ is symplectic
on all of~$P$, i.e., $\phi\in \operatorname{Hom}(P,P)$.  It remains to check that ${\textsf{c}}(\phi) =
\varphi$.  But the functor ${\textsf{c}}$ commutes with restrictions to
$P|_{U\cap {\mathaccent23{W}}}$ and
\begin{gather*}
{\textsf{c}}(\phi)|_{P|_{U\cap {\mathaccent23{W}}}} = {\textsf{c}}({\mathaccent23{\phi}}) = \varphi|{P|_{U\cap {\mathaccent23{W}}}}
\end{gather*}
by construction.  Hence, by Lemma~\ref{lem4:1}, ${\textsf{c}}(\phi) = \varphi$.
\end{proof}

\begin{Lemma}\label{lem4:3}
  Suppose $U\subset W$ is an open subset with $H^2(U,{\mathbb Z}) =0$.  Then for
    any $P_1, P_2 \in {\textsf{STB}}_\psi (U)$ the map
\begin{gather*}
{\textsf{c}}= {\textsf{c}}_U \colon \  \operatorname{Hom}(P_1, P_2) \to \operatorname{Hom}({\textsf{c}}(P_1) , {\textsf{c}}(P_2)),
\end{gather*}
is a bijection.
\end{Lemma}
\begin{proof}
By Lemma~\ref{lem4:1}, the map ${\textsf{c}}$ is an injection.
So we only need to check that ${\textsf{c}}$ is onto.
  Let $\varphi \in
  \operatorname{Hom}({\textsf{c}}(P_1), {\textsf{c}}(P_2))$.  By Lemma~\ref{loc-iso-stb} there exists an
  isomorphism $\phi\colon P_1\to P_2$.  Then ${\textsf{c}}(\phi)^{-1} \circ \varphi \in
  \operatorname{Hom}({\textsf{c}}(P_1), {\textsf{c}}(P_1))$. By Lemma~\ref{lem4:2}
\begin{gather*}
{\textsf{c}}(\phi)^{-1} \circ \varphi = {\textsf{c}} (\nu)
\end{gather*}
for some $\nu\in \operatorname{Hom}(P_1,P_1)$.  Hence
\begin{gather*}
\varphi = {\textsf{c}}(\phi)\circ {\textsf{c}} (\nu) = {\textsf{c}}(\phi\circ \nu).\tag*{\qed}
\end{gather*}
\renewcommand{\qed}{}
\end{proof}

We are now in position to f\/inish the proof that ${\textsf{c}}$ is fully faithful
by observing that for any two objects $P_1, P_2\in {\textsf{STB}}_\psi (W)$ the
functions $\underline{\operatorname{Hom}} (P_1, P_2)$ and $\underline{\operatorname{Hom}}
({\textsf{c}}(P_1), {\textsf{c}}(P_2))$ from the collection of open subset of~$W$ to sets
given respectively by
\begin{gather*}
\underline{\operatorname{Hom}} (P_1, P_2)(U):= \operatorname{Hom} (P_1|_U, P_2|_U).
\end{gather*}
and
\begin{gather*}
\underline{\operatorname{Hom}} ({\textsf{c}}(P_1), {\textsf{c}}(P_2))(U):= \operatorname{Hom} ({\textsf{c}}(P_1)|_U, {\textsf{c}}(P_2)|_U)
\end{gather*}
are sheaves.  Moreover
\begin{gather*}
{\textsf{c}} = {\textsf{c}}_U \colon \  \operatorname{Hom} (P_1|_U, P_2|_U) \to \operatorname{Hom} ({\textsf{c}}(P_1)|_U, {\textsf{c}}(P_2)|_U)
\end{gather*}
is a map of sheaves.  By Lemma~\ref{lem4:3} the map $c_U$ is a
bijection for any contractible open set $U$.  Hence ${\textsf{c}}\colon
\underline{\operatorname{Hom}} (P_1, P_2) \to \underline{\operatorname{Hom}} ({\textsf{c}}(P_1), {\textsf{c}}(P_2))$
is an isomorphism of sheaves.

This proves that for {\em any} open subset $U\subset W$ the functor
\begin{gather*}
{\textsf{c}}_U\colon \  {\textsf{STB}}_\psi(U)\to {\textsf{STM}}_\psi (U)
\end{gather*}
is fully faithful.  It remains to prove that ${\textsf{c}}$ is essentially
surjective.  As a f\/irst step in the proof of essential surjectivity we
observe that the objects on ${\textsf{STB}}_\psi(W)$ satisfy descent in
the sense of Grothendieck:

\begin{Lemma}\label{lem4:5}\label{lem:descent_for_stb}\sloppy
  Let $\{U_i\}_{i\in I}$ be an open cover of the manifold with corners
  $W$, $U_{ij}:= U_i \cap U_j$ and $U_{ijk} = U_i \cap U_j\cap U_k$
  for all $i,j,k\in I$.  Suppose we have a collection of objects $P_i
  \in {\textsf{STB}}_\psi(U_i)$ and isomorphisms $\Phi_{ij}\colon P_j|_{U_{ij}}\to
  P_i|_{U_{ij}}$ defining a $($normalized$)$ cocycle: $\Phi_{ii} = {\rm id}$,
  $\Phi_{ji} = \Phi_{ij}^{-1} $ and
\begin{gather*}
\Phi_{ij}|_{U_{ijk}}\circ \Phi_{jk}|_{U_{ijk}}\circ \Phi_{ki}|_{U_{ijk}} = {\rm id}
\end{gather*}
for all triples $i,j,k\in I$.  Then there exists an object
$P\in {\textsf{STB}}_\psi (W)$ and isomorphisms \mbox{$\gamma_i\colon  P|_{U_i}\to P_i$} so that
\begin{gather}\label{eq:descent}
\xy
(-20,10)*+{P_j|_{U_{ij}} }="1";
(20,10)*+ {P|_{U_{ij}}}="2";
(-20,-7)*+{P_i|_{U_{ij}}}="3";
(20,-7)*+{P|_{U_{ij}}}="4";
 {\ar@{<-}^{\gamma_j} "1";"2"};
{\ar@{=}_{}  "2";"4"};
 {\ar@{->}_{\Phi_{ij}} "1";"3"};
{\ar@{<-}^{\gamma_i} "3";"4"};
\endxy
\end{gather}
commutes.
\end{Lemma}

\begin{proof}
  We may take $P= \big(\bigsqcup_{i\in I} P_i\big)/_{\sim}$
  where $\sim$ is the
  equivalence relation def\/ined by the $\Phi_{ij}$s.  Then $P$ is a
  principal $G$-bundle over $W$ and the symplectic $G$-invariant forms
  on the $P_i$s def\/ine a $G$-invariant symplectic form on $P$.  The maps
  $\gamma_i^{-1} \colon P_i \to P|_{U_i}$ are induced by the inclusions
  $P_i\hookrightarrow \bigsqcup_{j\in I} P_j$.
\end{proof}

\begin{Lemma}\label{lem4:6} For any open subset $U$ of $W$
the functor ${\textsf{c}}\colon {\textsf{STB}}_\psi (U)\to {\textsf{STM}}_\psi (U)$ is essentially surjective.
\end{Lemma}

\begin{proof}
  Given $M\in {\textsf{STM}}_\psi(U)$, we want to show that it is isomorphic to
  ${\textsf{c}}(P)$ for some $P\in {\textsf{STB}}_\psi (U)$.

  Since ${\textsf{STB}}_\psi (U)$ is nonempty, we may choose an object $P'\in
  {\textsf{STB}}_\psi(U)$.  By Lemma~\ref{lem:loc-iso} ${\textsf{c}}(P')$ and $M$ are
  locally isomorphic. Therefore there is a cover $\{U_i\}_{i\in I}$ of
  $U$ and a family of isomorphisms $\{\varphi_i\colon {\textsf{c}}(P')|_{U_i}\to
  M|_{U_i}\}$.  Set
\begin{gather*}
P_i := P'|_{U_i}.
\end{gather*}
Consider the collection of isomorphisms
\begin{gather*}
\varphi_{ij} := (\varphi_i|_{U_{ij}})^{-1} \circ \varphi_j|_{U_{ij}}\colon \
{\textsf{c}}(P_j)|_{U_{ij}}\to {\textsf{c}}(P_i)|_{U_{ij}}, \qquad i,j\in I.
\end{gather*}
Since ${\textsf{c}}$ is fully faithful, there are unique isomorphisms
\begin{gather*}
\Phi_{ij}\colon \ P_j|_{U_{ij}}\to P_i|_{U_{ij}}
\end{gather*}
with ${\textsf{c}}(\Phi_{ij})= \varphi_{ij}$.  Since ${\textsf{c}}$ commutes with
restrictions to open subsets and since $\{\varphi_{ij}\}_{i,j\in I}$
form a cocycle and the $\Phi_{ij}$ are unique,
$\{\Phi_{ij}\}_{i,j\in I}$ form a cocycle as well.  By
Lemma~\ref{lem4:5} there is  $P\in {\textsf{STB}}_\psi(W)$ and a family
$\{\gamma_i\colon  P|_{U_i} \to P_i\}$ of isomorphisms so that~\eqref{eq:descent} commutes.  Then
\begin{gather*}
\xy
(-60,10)*+{M|_{U_{ij}} }="a";
(-60,-7)*+{M|_{U_{ij}}}="b";
(-20,10)*+{{\textsf{c}}(P_j)|_{U_{ij}} }="1";
(20,10)*+ {{\textsf{c}}(P)|_{U_{ij}}}="2";
(-20,-7)*+{{\textsf{c}}(P_i)|_{U_{ij}}}="3";
(20,-7)*+{{\textsf{c}}(P)|_{U_{ij}}}="4";
 {\ar@{<-}^{{\textsf{c}}(\gamma_j)} "1";"2"};
{\ar@{=}_{}  "2";"4"};
 {\ar@{->}_{\varphi_{ij}} "1";"3"};
{\ar@{<-}^{{\textsf{c}}(\gamma_i)} "3";"4"};
{\ar@{=}_{}  "a";"b"};
 {\ar@{->}_{\varphi_{j}} "1";"a"};
 {\ar@{->}_{\varphi_{i}} "3";"b"};
\endxy
\end{gather*}
commutes as well.  Consequently
\begin{gather*}
\varphi_i \circ {\textsf{c}}(\gamma_i)|_{{\textsf{c}}(P)|_{U_{ij}}} = \varphi_j \circ
{\textsf{c}}(\gamma_j)|_{{\textsf{c}}(P)|_{U_{ij}}}.
\end{gather*}
Therefore the family $\{\varphi_i \circ {\textsf{c}}(\gamma_i)\colon
{\textsf{c}}(P)|_{U_{i}}\to M_{U_i}\}$ gives rise to a well def\/ined isomorphism
${\textsf{c}}(P)\to M$.
\end{proof}

This completes our proof of Theorem \ref{thm:iso}.  In fact, we have
proved more:

\begin{Theorem} Let $\psi\colon W\to {\mathfrak g}^*$ be a u.l.e.
Then the functor
\begin{gather*}
{\textsf{c}} \colon \ {\textsf{STB}}_\psi \to {\textsf{STM}}_\psi
\end{gather*}
is an isomorphism of stacks over the site of open subsets of the
manifold with corners~$W$.
\end{Theorem}

\begin{proof}
Recall that
${\textsf{c}} \colon {\textsf{STB}}_\psi \to {\textsf{STM}}_\psi$ commutes with restrictions,
hence, it is a map of stacks.
By Theorem~\ref{thm:iso}, for every open subset $U \subset W$,
the functor $ c_U\colon  {\textsf{STB}}_\psi(U)\to {\textsf{STM}}_\psi (U) $
is an equivalence of categories.  Hence
${\textsf{c}} \colon {\textsf{STB}}_\psi \to {\textsf{STM}}_\psi$ is an isomorphism of stacks.
\end{proof}

\section[Characteristic classes and classif\/ication
of symplectic toric $G$-manifolds]{Characteristic classes and classif\/ication\\
of symplectic toric $\boldsymbol{G}$-manifolds}
\label{sec:invariants}

As we have seen in the previous section the functor
${\textsf{c}}\colon {\textsf{STB}}_\psi(W)\to {\textsf{STM}}_\psi(W)$ is an equivalence of categories.
Hence it def\/ines a bijection $\pi_0({\textsf{c}})\colon \pi_0({\textsf{STB}}_\psi(W))\to
\pi_0({\textsf{STM}}_\psi(W))$ between the sets of equivalence classes.  Thus to
f\/inish the proof of Theorem~\ref{thm:main}(\ref{2}) it is enough to
construct a bijection
\begin{gather*}
\pi_0({\textsf{STB}}_\psi (W))\  \leftrightarrow \
H^2(W, {\mathbb R})\times H^2 (W, {\mathbb Z}_G).
\end{gather*}
Recall that for a torus $G$ with integral lattice ${\mathbb Z}_G$ and a~manifold with corners $N$ there is a~bijection
\begin{gather*}
c_1\colon \  \pi_0 (BG(N))\to H^2 (N, {\mathbb Z}_G),
\end{gather*}
where $BG(N)$ denotes the category of principal $G$-bundles over $N$
and $c_1$ assigns to each isomorphism class $[P]\in \pi_0 (BG(N))$ of
a bundle $P$ its f\/irst Chern class $c_1(P)$.  Recall also that the
map $c_1$ is an isomorphism of presheaves.  Namely if
$V\stackrel{i}{\hookrightarrow} U\hookrightarrow N$ are two open
subsets of $N$ then the diagram
\begin{gather*}
\xy
(-20,10)*+{\pi_0(BG(U)) }="1";
(20,10)*+ {H^2(U,{\mathbb Z}_G)}="2";
(-20,-7)*+{\pi_0(BG(V)) }="3";
(20,-7)*+{H^2(V,{\mathbb Z}_G)}="4";
 {\ar@{->}^{c_1} "1";"2"};
{\ar@{->}_{i^*}  "2";"4"};
 {\ar@{->}_{i^*} "1";"3"};
{\ar@{->}^{c_1} "3";"4"};
\endxy
\end{gather*}
commutes.
Pre-composing with the map
$ \pi_0 ({\textsf{STB}}_\psi (\cdot)) \to \pi_0 (BG(\cdot))$
that is induced by the forgetful functor,
we get a homomorphism
\begin{gather} \label{c1}
c_1 \colon \  \pi_0 ({\textsf{STB}}_\psi(\cdot)) \to H^2(\cdot,{\mathbb Z}_G)
\end{gather}
of presheaves on $W$.

\begin{Proposition}\label{prop:char_class}
  Fix a u.l.e.\ $\psi\colon W\to {\mathfrak g}^*$.  The homomorphism~\eqref{c1}
  extends to an isomorphism of presheaves
\begin{gather}\label{eq:c_1-c_hor}
  (c_1, c_{\text{\rm hor}})\colon \ \pi_0 ({\textsf{STB}}_\psi(\cdot))\to
H^2(\,\cdot\,, {\mathbb Z}_G)\times H^2 (\,\cdot\,, {\mathbb R}).
\end{gather}
\end{Proposition}

\begin{Definition}\sloppy
  We call the second component of the isomorphism~\eqref{eq:c_1-c_hor}
  the {\em horizontal class}. We say informally that
  $c_{\text{hor}}([(P,\omega)])$ is the horizontal class of the symplectic toric
bundle \mbox{$(P\to W, \omega)$}.
\end{Definition}

\begin{proof}[Proof of \protect{Proposition~\ref{prop:char_class}}]
Fix an open set $U\subset W$.  We construct a bijection
\begin{gather*}
F=F_U\colon \  H^2 (U, {\mathbb Z}_G)\times H^2 (U,{\mathbb R})\to
\pi_0({\textsf{STB}}_\psi (U))
\end{gather*}
that commutes with pullbacks by open inclusions $i\colon V\hookrightarrow U$
and which is the inverse of the map~\eqref{eq:c_1-c_hor} on~$U$.
Here, we take $H^2(\cdot,{\mathbb R})$ to be the second de Rham cohomology.

Given $(c, [\beta]) \in H^2 (U, {\mathbb Z}_G)\times H^2 (U,{\mathbb R}) $ choose a
principal $G$-bundle $P$ with $c_1 (P) = c$.  By
Lemma~\ref{lemma:nondeg} a choice of a connection 1-form $A\in
\Omega^1(P,{\mathfrak g})^G$ def\/ines a symplectic form
\begin{gather*}
\omega_{A,\beta} = d \langle \mu, A\rangle + \pi^*\beta.
\end{gather*}
If $\beta' \in [\beta]$ is another closed 2-form representing the
class $[\beta]\in H^2(U, {\mathbb R})$ then $\beta' = \beta + d\gamma$ for some
$\gamma\in \Omega^1(U)$.  By Lemma~\ref{unique omega} the objects
$(P,\omega_{A,\beta})$ and $(P,\omega_{A,\beta'})$ are isomorphic in
${\textsf{STB}}_\psi (U)$.  If $A'$ is a dif\/ferent choice of a connection on $P$
then $A-A' = \pi^*a$ for some $a\in \Omega^1 (U,{\mathfrak g})$.  Consequently
by Lemma~\ref{unique omega} the symplectic bundles  $(P,\omega_{A, \beta})$ and
$(P,\omega_{A',\beta})$ are also isomorphic.  We conclude that the map
$F\colon H^2 (U, {\mathbb Z}_G)\times H^2 (U,{\mathbb R}) \to \pi_0({\textsf{STB}}_\psi (U))$ that
assigns to a pair $(c,[\beta])$ the isomorphism class of
$(P,\omega_{A,\beta})$ with ${\textsf{c}}_1(P) = c$ is well-def\/ined.

Given any $(P,\omega)\in {\textsf{STB}}_\psi(U)$, Lemma~\ref{unique omega}
implies that $\omega= \omega_{A,\beta}$ for some closed 2-form $\beta$
on~$U$.  Hence $F$ is onto.

Suppose $F(c, [\beta]) = [(P,\omega)]= F(c',[\beta'])$ for some
$(c, [\beta]), (c',[\beta'])\in  H^2 (U, {\mathbb Z}_G)\times H^2 (U,{\mathbb R})$.  Then
\begin{gather*}
c = c_1(P) = c'.
\end{gather*}
Now $F(c, [\beta]) =[(P, \omega_{A,\beta})]$ and $F(c', [\beta'])
=[(P, \omega_{A',\beta'})]$ for some connections
$A,A'\in\Omega^1(P,{\mathfrak g})$.  Since $[(P, \omega_{A,\beta})]= [(P,
\omega_{A',\beta'})]$ there is a gauge transformation $f\colon P\to P$ with
$f^*\omega _{A',\beta'} = \omega_{A,\beta}$.  Since $A$ and $f^*A'$ are
  both connections on $P$,
\begin{gather*}
f^*A' - A = \pi^* a
\end{gather*}
for some $a\in \Omega^1 (U,{\mathfrak g})$.  Since $\mu \circ f = \mu$ and $\pi
\circ f = \pi$,
\begin{gather*}
f^*(\omega_{A',\beta'})= f^*( d \langle \mu , A'\rangle +\pi^*\beta')
= d \langle \mu\circ f , f^*A'\rangle +f^*\pi^*\beta' = d \langle \mu,
f^*A'\rangle +\pi^*\beta'.
\end{gather*}
Consequently,
\begin{gather*}
0 = f^*(\omega_{A',\beta'}) - \omega_{A,\beta} = d\langle \mu,
f^*A'\rangle +\pi^*\beta' - d\langle \mu, A\rangle -\pi^*\beta.
\end{gather*}
Therefore
\begin{gather*}
\pi^*(\beta -\beta') =d\langle \mu,
f^*A' -A \rangle = d \langle \mu , \pi^*a\rangle
= \pi^* (d\langle \psi, a\rangle).
\end{gather*}
Hence
\begin{gather*}
\beta -\beta' = d\langle \psi, a\rangle.
\end{gather*}
Therefore $[\beta] = [\beta']$, and $F$ is one-to-one.
\end{proof}

As an immediate consequence we have
\begin{proof}[Proof of \protect{Theorem~\ref{thm:main}(\ref{2})}]
It follows from Proposition~\ref{prop:char_class} that the composite
\begin{gather*}
 \pi_0({\textsf{STM}}_\psi(W))\xrightarrow{ \pi_0({\textsf{c}})^{-1} } \pi_0({\textsf{STB}}_\psi (W))
 \xrightarrow{(c_1, c_{\text{hor}})} H^2 (W, {\mathbb Z}_G)\times H^2 (W,{\mathbb R})
\end{gather*}
is a bijection.
\end{proof}

\begin{Definition}[Chern and horizontal classes of a symplectic toric manifold]
  Let $(\pi\colon M\to W, \omega) \in {\textsf{STB}}_\psi(W)$ be a symplectic toric
  manifold over a u.l.e.\ $\psi\colon W\to {\mathfrak g}^*$. We def\/ine its {\em Chern
    class} to be the f\/irst Chern class of the corresponding principal
  torus bundle:
\begin{gather*}
c_1(M, \omega,\pi) := c_1 \circ  \pi_0({\textsf{c}})^{-1} ([M,\omega,\pi]).
\end{gather*}
Similarly its {\em horizontal class} is the horizontal class of the
corresponding bundle:
\begin{gather*}
c_{\text{hor}}(M, \omega,\pi) := c_{\text{hor}}\circ  \pi_0({\textsf{c}})^{-1} ([M,\omega,\pi]).
\end{gather*}
\end{Definition}
Another consequence of Proposition~\ref{prop:char_class} is
\begin{Corollary}\label{STB_is_a_gerbe}
  Fix a u.l.e.\ $\psi\colon W\to {\mathfrak g}^*$. If $H^2(W,{\mathbb Z})= 0$ then any two
  objects of ${\textsf{STM}}_\psi(W)$ are isomorphic.
\end{Corollary}
This corollary signif\/icantly strengthens Lemma~\ref{lem:loc-iso}.

\begin{Remark} The Chern and horizontal classes for a symplectic toric
  manifold over $W$ have a~nice geometric interpretation in terms of
  the restriction to the interior~${\mathaccent23{W}}$.

  If $W$ is a manifold with corners and ${\mathaccent23{W}}$ is its interior then
  the inclusion map $\iota\colon  {\mathaccent23{W}} \hookrightarrow W$ is a~smooth
  homotopy equivalence.  (A homotopy inverse is obtained from the f\/low
  of a vector f\/ield on~$W$ that is supported in a small neighborhood
  of the topological boundary~$\partial W$, is transverse to all the
  codimension 1 strata of~$W$, and points inward along the boundary~$\partial W$.)  So the restriction maps $\iota^*\colon  H^2(W;{\mathbb Z}_G) \to
  H^2({\mathaccent23{W}};{\mathbb Z}_G)$ and $\iota^*\colon H^2(W;{\mathbb R}) \to H^2({\mathaccent23{W}};{\mathbb R})$ are
  isomorphisms.  By Proposition~\ref{prop:char_class} the diagram
\begin{gather*}
\xy
(-20,10)*+{\pi_0 ({\textsf{STM}}_\psi (W)) }="1";
(40,10)*+ {H^2(W,Z_G)\times H^2(W, {\mathbb R})}="2";
(-20,-7)*+{\pi_0 ({\textsf{STM}}_\psi ({\mathaccent23{W}}))}="3";
(40,-7)*+{H^2({\mathaccent23{W}},Z_G)\times H^2({\mathaccent23{W}}, {\mathbb R})}="4";
 {\ar@{->}^{(c_1, c_{\text{hor}})} "1";"2"};
{\ar@{->}^{\iota^*}  "2";"4"};
 {\ar@{->}_{\iota^*} "1";"3"};
{\ar@{->}^{(c_1, c_{\text{hor}})} "3";"4"};
\endxy
\end{gather*}
commutes.  Since $c\colon {\textsf{STB}}_\psi ({\mathaccent23{W}}) \to {\textsf{STM}}_\psi({\mathaccent23{W}})$ is isomorphic
to the identity, the induced map~$\pi_0(c)$ is the identity.  It
follows that for a symplectic toric manifold $(M,\omega, \pi)\in
{\textsf{STM}}_\psi(W)$
\begin{gather*}
\iota^*c_1 (M,\omega, \pi) = c_1 (M|_{{\mathaccent23{W}}}, \omega, \pi) = c_1 (M|_{\mathaccent23{W}}),
\end{gather*}
where \looseness=1 $c_1 (M|_{\mathaccent23{W}})$ is the Chern class of the principal
$G$-bundle $M|_{\mathaccent23{W}} \to {\mathaccent23{W}}$. Note that this in particular relates
our def\/inition of the Chern class of a symplectic toric manifold to
Duistermaat's def\/inition of the Chern class of a completely integrable
system \cite{Duis} (cf.\ Remarks~\ref{rmrk:1.9} and~\ref{rmrk:1.11}).

\sloppy Similarly
\begin{gather*}
\iota^*c_{\text{hor}} (M,\omega, \pi) = c_{\text{hor}} (M|_{{\mathaccent23{W}}}, \omega, \pi),
\end{gather*}
where the class on the right is the horizontal class of the symplectic
principal $G$-bundle \mbox{$M|_{\mathaccent23{W}} \to {\mathaccent23{W}}$}.
\end{Remark}

\section{Toric manifolds determined by their moment map images}
\label{sec:the_rest}

As we mentioned in the introduction, in general the image of the
moment map doesn't tell us much about the symplectic toric manifold.
There are two reasons for this. First of all, the orbital moment map
may not be an embedding~-- see Example~\ref{ex:mmap-not-emb}.
Secondly, even when the orbital moment map is an embedding the second
integral cohomology of the orbit space, which then has to be
(isomorphic to) the image of the moment map, may not be trivial.

\begin{Example} Consider a three-dimensional torus $G$.  Then $H^2({\mathfrak g}^*
\setminus \{0\}, {\mathbb Z}) = {\mathbb Z}$.  By Theorem~\ref{thm:main}
there are ${\mathbb Z}_G \times {\mathbb R}$ isomorphism classes of symplectic toric
manifolds over ${\mathfrak g}^* \setminus \{0\}$ (they are all principal
$G$-bundles over ${\mathfrak g}^* \setminus \{0\}$).  For all of these manifolds
the orbital moment map is an embedding and the moment map image is
${\mathfrak g}^*\setminus \{0\}$. This family of manifolds has been studied by
Bates~\cite{Bates}.  He proves that any principal $G$-bundle
$P\to {\mathfrak g}^* \setminus \{0\}$ admits a symplectic form making the f\/ibers of
$P\to {\mathfrak g}^* \setminus \{0\}$ Lagrangian and the action of~$G$ on~$P$ Hamiltonian.
\end{Example}

\begin{Example}
  It is also easy to construct examples of symplectic
  toric manifolds where the orbital moment map is an embedding, the
  torus action is not free and the second integral cohomology of the
  orbit space is nontrivial.  For instance let $G$ again be a three-dimensional torus, let $\Delta \subset {\mathfrak g}^*$ be a unimodular
  simplex (or any other unimodular polytope) and let $W = \Delta
  \setminus \{w_0\}$, where $w_0$ is a point in the interior of
  $\Delta$.  Then again $H^2(W \setminus \{ w_0 \}, {\mathbb Z}) = {\mathbb Z}$ and
  consequently there are ${\mathbb Z}_G \times {\mathbb R}$ isomorphism classes of
  symplectic toric manifolds over $W$.  For every symplectic toric
  manifold $M$ over $W\hookrightarrow {\mathfrak g}^*$ the points over the
  vertices of $W$ are the f\/ixed points of the $G$ action.
\end{Example}

Theorem~\ref{thm:main} implies that (the isomorphism
class of) a symplectic toric manifold $(M,\omega, \mu)$ is uniquely
determined by $\mu(M)$ if (1) the orbital moment map is an
embedding and (2) $H^2(\mu(M),{\mathbb Z}) =0$.  The connectedness and
convexity theorems of Atiyah, Guillemin and Sternberg imply that any
compact connected symplectic toric manifold falls into this class of
toric manifolds. But there is more.  Recall that the assumption of
compactness in connectedness and convexity theorems can be weakened.
Namely,

\begin{Theorem}[cf.\ \protect{\cite[Theorem~4.3]{LMTW}}]\label{thm:LMTW}
    Let $\mu\colon M\to {\mathfrak g}^*$ be a moment map for an action of a~torus~$G$
    on a symplectic manifold $(M,\omega)$.  Suppose there exists a
    convex open subset $U\subset {\mathfrak g}^*$ such that $\mu(M)\subset U$ and
    $\mu\colon M\to U$ is proper.  Then the image $\mu(M)$ is convex and
    each fiber of $\mu$ is connected.
\end{Theorem}

It will be convenient to have the following def\/inition:

\begin{Definition}
  A map $F\colon M\to V$ from a space $M$ to a f\/inite-dimensional vector
  space $V$ is {\em proper as a map into a~convex open set} if
  there is a convex open set $U\subset V$ so that $F(M)\subset U$ and
  $F\colon M\to U$ is proper.
\end{Definition}

\begin{Proposition} \label{Delzant over subset}
Let $(M,\omega,\mu)$ be a connected symplectic toric $G$-manifold
whose moment map~$\mu$ is proper as a~map into a~convex open set.
Then the image~$\mu(M)$ determines~$(M,\omega,\mu)$ up to isomorphism.
\end{Proposition}

\begin{proof}
  Theorem~\ref{thm:LMTW} and Lemma~\ref{compactly generated}
  imply that the orbital moment map $\overline{\mu}$ is an embedding.
  Consequently, $\bar{\mu}\colon M/G\to \mu(M)$ is a dif\/feomorphism
  of manifolds with corners and \mbox{$\mu\colon M\to \mu(M)$} is a quotient map.
  That is $(M,\omega, \mu)$ is a symplectic toric manifold over
  $\mu(M)$.  By Theorem~\ref{thm:LMTW}, $\mu(M) \simeq M/G$ is
  convex, hence $H^2(M/G,{\mathbb Z}) = 0$.  It now follows from
  Theorem~\ref{thm:main} that any two symplectic toric manifolds over
  $\mu(M)$ are isomorphic.
\end{proof}

\begin{Remark}
  Proposition~\ref{Delzant over subset} can be used for
  ``coordinatization'' of compact symplectic toric manifolds in the
  sense of Duistermaat and Pelayo~\cite{DP}. Namely, let
  $(M,\omega,\mu)$ be a compact connected symplectic toric
  $G$-manifold with momentum polytope $\Delta = \mu(M)$.  Express
  $\Delta$ as the intersection of half-spaces $H_1,\ldots,H_N$ whose
  boundaries are the af\/f\/ine spans of the facets of $\Delta$.  For each
  vertex $\epsilon$ of $\Delta$, let $C_\epsilon$ be the intersection of those
  $H_j$ whose boundary contains $\epsilon$ (this is the tangent cone to
  $\Delta$ at $\epsilon$), and let ${\mathcal T}_\epsilon$ be the intersection of
  $\text{interior}(H_j)$ over those $j$ such that $\epsilon \in
  \text{interior}(H_j)$.

Let $U_\epsilon$ denote the preimage in $M$ of the convex open subset
${\mathcal T}_\epsilon$ of ${\mathfrak g}^*$.
The sets $U_\epsilon$ form a~covering of $M$ by $G$-invariant open dense subsets.
By Proposition~\ref{Delzant over subset},
each of these sub\-sets~$U_\epsilon$ is equivariantly symplectomorphic
to a $G$-invariant open subset of~${\mathbb C}^n$,
where $G$ acts on~${\mathbb C}^n$ through the isomorphism~$G \to (S^1)^n$
for which the momentum map image is~$C_\epsilon$.
\end{Remark}

The reader familiar with Delzant's paper may wonder which of the
symplectic toric manifolds that we classify can be obtained as
symplectic quotients of some standard ${\mathbb C}^N$ by an action of a
subtorus of the standard torus ${\mathbb T}^N = (S^1)^N$.
The following theorem and its proof are the result of our discussion
with Chris Woodward.  We thank Chris for bringing up the question
and helping us prove the answer.

\begin{Theorem} \label{reduction of CN}
  A symplectic toric $G$-manifold $(M,\omega,\mu \colon M\to {\mathfrak g}^*)$ is
  isomorphic to a regular symplectic quotient of ${\mathbb C}^N$ by a subtorus of the
  standard torus ${\mathbb T}^N$ if and only if
  its orbital moment map $\overline{\mu} \colon M/G \to {\mathfrak g}^*$
  is an embedding and its image is a closed convex
  polyhedral subset of~${\mathfrak g}^*$ with at least one vertex
  and at most~$N$  facets.
\end{Theorem}

\begin{Remark}
We already know that the orbital momentum map is locally
an embedding as a~manifold with corners.
So it is a~global embedding as a~manifold with corners
if and only if it is a~global embedding topologically.
\end{Remark}

\begin{proof}[Proof of Theorem~\ref{reduction of CN}]
  We f\/irst argue that the conditions of the theorem are necessary: if
  a~symplectic toric manifold $(M,\omega, \mu)$ is a non-singular
  symplectic quotient of ${\mathbb C}^N$ by a subtorus $K\hookrightarrow {\mathbb T}^N$,
  then the orbital moment map $\overline{\mu}$ is an embedding and the
  moment map image~$\mu(M)$ is a unimodular polyhedral subset of~${\mathfrak g}^*$ with at least one vertex.

  Recall that for the standard action of ${\mathbb T}^N$ on ${\mathbb C}^N$ the map
  $\psi \colon {\mathbb C}^N \to {\mathbb R}^N$,   $\psi(z_1, \ldots, z_N) = (|z_1|^2,
  \ldots, |z_N|^2)$, is a moment map (with an appropriate
  identif\/ication of the dual of the Lie algebra of ${\mathbb T}^N$ with
  ${\mathbb R}^N$).
  The image of $\psi$ is the positive orthant ${\mathbb R}^N_+ :=
  \{t\in {\mathbb R}^N\,|\, t_i\geq 0 \, \forall\, i\}$,
  and the orbital moment map
  $\overline{\psi} \colon {\mathbb C}^N/{\mathbb T}^N = {\mathbb R}^N_+\to {\mathbb R}^N$ is an embedding.
  The inclusion $K\hookrightarrow {\mathbb T}^N$ induces an
  inclusion $i^T \colon {\mathfrak k}\hookrightarrow \operatorname{Lie}({\mathbb T}^N)$ of Lie algebras
  and, dually, the projection $i^T \colon \operatorname{Lie}({\mathbb R}^N)^*= {\mathbb R}^N\to
  {\mathfrak k}^*$.  Then $\varphi:=i^T\circ \psi \colon {\mathbb C}^N\to {\mathfrak k}^*$ is a~moment map for the action of~$K$ on~${\mathbb C}^N$.  Suppose $\nu\in {\mathfrak k}^*$
  is a point such that the action of~$K$ on its preimage~$\varphi^{-1}
  (\nu)$ is free. Then~$\nu$ is necessarily a regular value of
  $\varphi$. Moreover, the symplectic quotient $M:= \varphi ^{-1}
  (\nu)/K$ is a~symplectic toric $G$-manifold for $G= {\mathbb T}^N/K$ (cf.~\cite{De}).  The moment map $\mu \colon M\to {\mathfrak g}^*$ has the
  following description.  Since ${\mathfrak g}^*$ is canonically isomorphic to
  the annihilator ${\mathfrak k}^\circ$ of ${\mathfrak k}$ in $\operatorname{Lie}({\mathbb T}^N)^* ={\mathbb R}^N$ and
  since $(i^T)^{-1} (\nu) = \lambda + {\mathfrak k}^\circ$ for some $\lambda \in
  (i^T)^{-1} (\nu)$, we can identify ${\mathfrak g}^*$ with the af\/f\/ine plane
  $(i^T)^{-1} (\nu)$.  With this identif\/ication, the restriction
\begin{gather*}
\psi|_{\varphi^{-1} (\nu)} = \psi|_{\psi^{-1} ((i^T)^{-1} (\nu))} \colon
\  \psi^{-1} \big(\big(i^T\big)^{-1} (\nu)\big)   \to   \big(i^T\big)^{-1} (\nu)
\end{gather*}
descends to $\mu \colon M=\varphi^{-1} (\nu)/K \to (i^T)^{-1} (\nu)\simeq
{\mathfrak g}^*$.  Since all the f\/ibers of $\psi$ are ${\mathbb T}^N$-orbits, the f\/ibers
of $\mu$ are ${\mathbb T}^N/K= G$-orbits.  Since $\overline{\psi}$ is an open map
to its image, so is $\overline{\mu}$.
Hence the orbital moment map
$\overline{\mu} \colon M/G\to (i^T)^{-1} (\nu) \simeq {\mathfrak g}^*$ is an embedding.
Since the image of~$M$ under~$\mu$ is $\psi({\mathbb C}^N) \cap (i^T)^{-1} (\nu)
= {\mathbb R}^N_+ \cap (i^T)^{-1} (\nu)$, it has at most $N$ facets.  Now we
argue that the image has at least one vertex.  We give a~symplectic (as opposed to a convex geometry) argument.

The function $f(z) = \sum |z_j|^2 \colon {\mathbb C}^N\to {\mathbb R}$ is proper, non-negative,
${\mathbb T}^N$-invariant, and its Hamiltonian vector f\/ield generates a circle action.
Hence $f|_{\varphi^{-1} (\nu)}$ descends to
a proper non-negative $G$-invariant periodic Hamiltonian $\overline{f}$ on
$M$.  Since $\overline{f}$ is non-negative and proper, it achieves a~minimum somewhere on~$M$.  Since~$\overline{f}$ is periodic and
$G$-invariant, the set of points where~$\overline{f}$ achieves its minimum
is a $G$-invariant symplectic submanifold~$M'$ of~$M$.  Since $\overline{f}$
is proper, $M'$~is compact.  Since the action of~$G$ on~$M'$ is
Hamiltonian, it has to have a f\/ixed point.  The image of this f\/ixed
point under the moment map~$\mu$ is a vertex of the polyhedral set~$\mu(M)$.

Now suppose that the orbital moment map $\overline{\mu} \colon M/G \to {\mathfrak g}^*$ of a~symplectic toric $G$-manifold $(M,\omega,\mu \colon M\to {\mathfrak g}^*)$ is an
embedding, that its image $A =\overline{\mu}(M/G)$ is an intersection of
$N$ closed half-spaces, and that $A$ has a vertex $*$.
Let $H_1, \ldots, H_n$ denote the half-spaces
whose supporting hyperplanes form the facets meeting at~$*$.
Since $A$ is unimodular, $n=\dim G$
and the intersection $H_1\cap \cdots \cap H_n$ is isomorphic (as a~unimodular cone) to ${\mathbb R}^n_+$.  (In particular, $N \geq \dim G$.)
Consequently there is an isomorphism
$G\to {\mathbb T}^n$ so that the image of a corresponding moment map $\psi \colon {\mathbb C}^n
\to {\mathfrak g}^*$ is $H_1\cap \cdots \cap H_n$. Let $H_{n+1},\ldots, H_N$ denote
the remaining half-spaces, so that $A= H_1\cap\cdots\cap H_N$.  We now
successively apply the symplectic cut construction to ${\mathbb C}^n$ using the
half-spaces $H_{n+1},\ldots, H_N$~\cite{Lcuts}.  This amounts to
taking a symplectic quotient of ${\mathbb C}^n\times {\mathbb C}^{N-n}$ by an action of
$(S^1)^{N-n}$.  Since $A$ is unimodular, the quotient is smooth~\cite{LMTW}.  This symplectic quotient of ${\mathbb C}^N$ is
a symplectic toric $G$-manifold with moment map image $A$
and whose orbital moment map is an embedding.
Because~$A$ is convex, its second cohomology is trivial;
by Theorem~\ref{thm:main}, this symplectic quotient
of ${\mathbb C}^N$ is isomorphic to $(M,\omega,\mu \colon M \to {\mathfrak g}^*)$.
\end{proof}

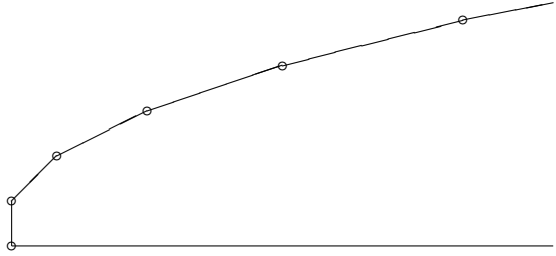
\begin{figure}[t]\centering
\setlength{\unitlength}{.6cm}
\begin{picture}(12,6)(0,0)
\put(0,0){\circle{.2}}
\put(0,0){\line(1,0){12}}
\put(0,0){\line(0,1){1}}
\put(0,1){\circle{.2}}
\put(0,1){\line(1,1){1}}
\put(1,2){\circle{.2}}
\put(1,2){\line(2,1){2}}
\put(3,3){\circle{.2}}
\put(3,3){\line(3,1){3}}
\put(6,4){\circle{.2}}
\put(6,4){\line(4,1){4}}
\put(10,5){\circle{.2}}
\put(10,5){\line(5,1){2}}
\end{picture}
\caption{A noncompact symplectic toric manifold
that is not a reduction of ${\mathbb C}^N$.}
\label{fig:example}
\end{figure}

We end this section with an example of a noncompact symplectic
toric $G$-manifold that is \emph{not} isomorphic to a symplectic
reduction of any ${\mathbb C}^N$:

\begin{Example}\label{ex:2.15}
Let $G = (S^1)^2$.
Let $W$ be the closed region in ${\mathfrak g}^* = {\mathbb R}^2$ that is bounded
on the bottom by the positive $x$ axis and on the top by the
polygonal path that is obtained by connecting, in this order,
the points  $(k(k-1)/2,k)$ for $k \in \{ 0,1,2,\ldots \}$.
See Fig.~\ref{fig:example}.
The set~$W$ is locally unimodular:
near the vertex $(0,0)$ it coincides with the positive orthant,
and near the vertex $v = (k(k-1)/2,k)$ for $k \geq 1$ it coincides with the cone
$v + {\mathbb R}_+ (-(k-1),-1) + {\mathbb R}_+ (k,1)$,
which is unimodular because $\det \left[ \begin{smallmatrix}
-(k-1) & -1 \\ k & 1 \end{smallmatrix} \right] = 1$.
By Theorem~\ref{thm:main}, there exists
a symplectic toric $G$-manifold $M$
with moment image~$W$.
Because~$W$ has inf\/initely many facets (edges)
and by Theorem~\ref{reduction of CN}, $M$ is not isomorphic
to a symplectic quotient of any~${\mathbb C}^N$.
\end{Example}

\appendix

\section{Manifolds with corners} \label{app:mfld-w-corner}

We quote Joyce \cite{Joyce}:
\begin{quote}
  ``[manifolds with corners] were f\/irst developed by Cerf~[1] and Douady~[2] in 1961, who were primarily interested in their
  Dif\/ferential Geometry. J\"anich [5] used manifolds with corners to
  classify actions of transformation groups on smooth manifolds.
  Melrose~[12, 13] and others study analysis of elliptic operators on
  manifolds with corners.
  \ldots\
  How one sets up
  the theory of manifolds with corners is not universally agreed, but
  depends on the applications one has in mind. \ldots
  there are at least four inequivalent def\/initions of manifolds with
  corners, two inequivalent def\/initions of boundary, and (including
  ours) four inequivalent def\/initions of smooth map in use in the
  literature.''
\end{quote}
The purpose of the appendix is to spell out our approach to manifolds
with corners and their maps.  In particular we spell out what we mean
for a~subset~$Y$ of a manifold with corners~$X$ to be naturally a
smooth manifold (Def\/inition~\ref{submanifold}) and what we mean by an
embedding of a~manifold with corners into a~manifold
(Def\/inition~\ref{def:embedding}).

\begin{Definition}[manifold with corners] \label{md w corners}\sloppy
Let $V$ be an (arbitrary) subset of ${\mathbb R}^n$.  A~map \mbox{$\varphi \colon V \to {\mathbb R}^m$}
is \emph{smooth} if for every point~$p$ of $V$ there exist
an open subset~$\Omega$ in~${\mathbb R}^n$ containing~$p$ and a~smooth map
from $\Omega$ to ${\mathbb R}^m$ whose restriction to~$\Omega \cap V$
coincides with $\varphi|_{\Omega \cap V}$.
A~map~$\varphi$ from~$V$ to a~subset of ${\mathbb R}^m$ is \emph{smooth}
if it is smooth as a~map to~${\mathbb R}^m$.  A~map~$\varphi$ from a~subset
of~${\mathbb R}^n$ to a~subset of~${\mathbb R}^m$ is a \emph{diffeomorphism}
if it is a bijection and both it and its inverse are smooth.

A \emph{sector} is the set $[0,\infty)^k \times {\mathbb R}^{n-k}$
where $n$ is a non-negative integer and $k$ is an integer between $0$ and $n$.
Let $X$ be a Hausdorf\/f second countable topological space.
A \emph{chart} on an open subset $U$ of $X$ is a homeomorphism $\varphi$
from $U$ to an open subset $V$ of a sector.
Charts $\varphi \colon U \to V$ and $\varphi' \colon U' \to V'$
are \emph{compatible} if $\varphi' \circ \varphi^{-1}$ is
a dif\/feomorphism from $\varphi(U \cap U')$ to $\varphi'(U \cap U')$.
An \emph{atlas} on $X$ is a set of pairwise compatible charts
whose domains cover $X$. Two atlases are \emph{equivalent}
if their union is an atlas.
A \emph{manifold with corners} is a Hausdorf\/f second countable
topological space equipped with an equivalence class of atlases.
\end{Definition}

We sometimes refer to an ordinary manifold as a ``manifold
without boundary or corners''.

\begin{Definition}[smooth map]
Let $X$ and $Y$ be manifolds with corners.
A map $h \colon X \to Y$ is \emph{smooth}
if for every point in $X$ there exists an open neighbourhood~$U$ in $X$
and an open subset~$U'$ of~$Y$
and charts $\varphi \colon U \to V$ and $\varphi' \colon U' \to V'$
of~$X$ and~ $Y$
such that~$h(U) \subset U'$
and such that $\varphi' \circ h \circ \varphi^{-1} \colon V \to V'$ is smooth.

A map $f$ from an (arbitrary) subset $A \subset X$ to $Y$ is smooth
if for every point in $A$ there exists a neighbourhood ${\mathcal O}$ in $X$
and a smooth map to $Y$ whose restriction to $A \cap {\mathcal O}$ coincides
with $f$; a map from $A$ to a subset $B$ of $Y$ is smooth if it is smooth
as a map to $Y$, and it is a~\emph{diffeomorphism} if it is smooth and has
a smooth inverse.
\end{Definition}

It is easy to check that the composition of two smooth maps is again smooth.
Hence manifolds with corners form a category.
The isomorphisms in this category are \emph{diffeomorphisms}.
We may refer to a smooth map between two manifolds with corners
as a {\em map of manifolds with corners}.

\begin{Definition} \label{strata} The \emph{dimension} of a manifold
  with corners $X$ is $n$ if the charts take values in sectors in
  ${\mathbb R}^n$.  A point $x$ of $X$ has \emph{index~$k$} if there exists a~chart $\varphi$ from a neighbourhood of~$x$ to $[0,\infty)^k \times
  {\mathbb R}^{n-k}$ such that $\varphi(x) = 0$; the index of a point is well
  def\/ined.  The \emph{$k$-boundary},~$X^{(k)}$, of $X$ is the set of
  points of index $\geq k$.  The (topological) \emph{boundary} of~$X$ is the
  $1$-boundary, $\partial X := X^{(1)}$.  The \emph{interior} of~$X$ is
  the complement of the boundary: ${\mathaccent23{X}} := X \setminus \partial X$; it is
  the set of points of index~0.
We refer to the connected components of the sets
\begin{gather*}
\partial^{(k)}X:= X^{(k)} \setminus X^{(k+1)}
\end{gather*}
 as the \emph{strata} of $X$.
\end{Definition}

\begin{Definition}\label{differential forms}
The \emph{tangent space} $T_xX$ of a manifold with corners $X$ at a
point $x\in X$ is the space of derivations at $x$ of germs at $x$ of
smooth functions def\/ined near $x$.  Thus, the tangent space is a
vector space even if the point $x$ is in the boundary of $X$.

Similarly the tangent bundle $TX$ of a manifold with corners $X$ is a
vector bundle over $X$ as is the cotangent bundle $T^*X$ and its
exterior powers. The total spaces of $TX$ and $T^*X$ are manifolds
with corners (cf.~\cite[p.~19]{Michor}).

A {\em differential $k$-form} on a manifold with corners $X$
is a smooth section of the $k$th exterior power of its cotangent bundle.

Exterior derivative $d$ makes sense on manifolds with corners.  So do
closed forms and symplectic forms.
Thus, the usual notions of a Hamiltonian action and a moment map
extend to manifolds with corners without change.

We note that, on an open subset of a sector,
a closed form locally extends as a \emph{closed} form, because
it locally has a primitive and the primitive has a local smooth extension.
\end{Definition}

\begin{Definition}\label{def:embedding}
  A smooth map $f\colon N\to M$ of manifolds with corners is an {\em
    embedding} if it is a~topological embedding and the dif\/ferential
  $df_x\colon T_xM\to T_{f(x)}N$ of $f$ is injective at every point $x\in M$.
  Equivalently, $f$ is an embedding if $f \colon N \to f(N)$
  is a dif\/feomorphism.
\end{Definition}

\begin{Example}
The inclusion $[0,\infty)^k \hookrightarrow {\mathbb R}^k$ is an embedding.
\end{Example}

One can prove (q.v.\ \cite[p.~21]{Michor}):
\begin{Lemma}\label{lem:domain}
  A manifold with corners $M$ can be embedded in a manifold
  $\tilde{M}$ $($without corners$)$ of the same dimension.
\end{Lemma}

\begin{Definition}[domain]\label{def:domain}
  If a manifold with corners $M$ is embedded in a manifold~$\tilde{M}$
  (without corners) and if $\dim M = \dim \tilde{M}$, we say that $M$ is
  a~{\em domain} in~$\tilde{M}$.
\end{Definition}

\begin{Remark} \label{rmrk:A9}
If $M\subset \tilde{M}$ is a domain and $f\colon M\to V$ is a~smooth map to some f\/inite-dimensional vector space~$V$, then~$f$
  extends to a smooth map $\tilde{f}$ from some open neighbourhood of~$M$ in $\tilde{M}$.
  We refer to
  $\tilde{f}$ as an {\em extension} of $f$ to $\tilde{M}$.
\end{Remark}

It is not true that the closure of a stratum of a manifold with corners
is a manifold with corners.  See for example Fig.~2.1 in~\cite{Joyce}.  Nor is it true a stratum
of codimension~$k$ lies in the closure of exactly~$k$ codimension~1
strata.

\begin{Definition} \label{def:faces} A manifold with corners~$X$ is a
{\em manifold with faces} (q.v.~\cite{Janich}) if every point of~$X$
of index~$k$ lies in the closure of exactly~$k$ codimension~1 strata.

We refer to the closures of the strata of a manifold with faces~$X$
as {\em faces} and
the codimension~1 faces as {\em facets.}
\end{Definition}

One can show that for a manifold with faces the
closure of a stratum is a manifold with faces (op.\ cit.).

\begin{Example} A unimodular cone is a manifold with faces.
\end{Example}

\begin{Example} A sector $[0,\infty)^k \times {\mathbb R}^{n-k}$ is a manifold
with faces.
\end{Example}

\begin{Remark}
  An open subspace of a manifold with faces is again a manifold with
  faces.  Consequently any manifold with corners is a manifold with
  faces locally: for any manifold with corners~$N$ and any point $x\in
  N$ there is an open neighbourhood~$U$ of~$x$ in~$N$ such that~$U$ is
  a~manifold with faces.  We will refer to such neighbourhood~$U$ as a~{\em neighbourhood with faces.}
\end{Remark}

The following def\/inition is nonstandard but is essential for the
purposes of this paper.

\begin{Definition} \label{submanifold}
  We say that a subset $Y$ of a manifold with corners $X$
  is \emph{naturally a smooth manifold}
  if it has a manifold structure such that the inclusion map
  $Y \hookrightarrow X$ is an embedding
  in the sense of Def\/inition~\ref{def:embedding}.
  If such a manifold structure on $Y$ exists, then it is unique.
\end{Definition}

\begin{Example}
  With Def\/inition~\ref{submanifold}, the parabola $\{(x,y)\in {\mathbb R}^2
  \,|\, y = x^2 \}$, as a subset of the upper half plane $\{(x,y)\in
  {\mathbb R}^2 \,|\, y \geq 0 \}$, is naturally a smooth manifold.
\end{Example}

\begin{Definition} \label{quotient}
  An \emph{action} of a Lie group $G$ on a manifold with corners $X$
  is a~homomor\-phism~$\rho$ from $G$ to the group of dif\/feomorphisms of
  $X$ such that the map $G \times X \to X$ given by $(a,x) \mapsto
  \rho(a)(x)$ is smooth.

  Given an action of a compact Lie group $G$ on a manifold with corners~$X$,
  we say that a~smooth
  map $\pi \colon X \to W$ from $X$ to another manifold with corners~$W$ is \emph{a quotient map} if for every $G$-invariant smooth map
  $f \colon X \to Y$ there exists a unique smooth map $\overline{f} \colon W
  \to Y$ such that $f = \overline{f} \circ \pi$.
  Such a map $\pi$ identif\/ies $W$ with the set $X/G$ of $G$ orbits.

  We say that $X/G$ is \emph{naturally a manifold with corners}
  if it has a manifold with corners structure such that the map $X \to X/G$
  that takes each point to its orbit is a quotient map in the above sense.
  If such a structure exists then it is unique.  We say that $X/G$
  is \emph{naturally a~smooth manifold}, or, simply, that it is a manifold,
  if it has such a structure without boundary or corners.
\end{Definition}

\begin{Definition}\label{def:pric_G-bundle_corners}
  Let $W$ be a manifold with corners and $G$ a Lie group.
  A \emph{principal $G$-bundle} over $W$ is a manifold with corners $P$
  equipped with a right action of $G$ and with a map $\pi\colon P\to W$
  making it into a topological principal $G$-bundle
  in which local trivializations can be chosen to be equivariant
  dif\/feomorphisms of manifolds with corners.
\end{Definition}

On a manifold with corners~$M$, the de Rham cohomology
is well def\/ined and invariant under homotopy.
The Poincar\'e lemma holds: every closed form is locally exact.
The proofs are exactly as for ordinary manifolds.
A key step is that if $\beta$ is a $k$-form on $[0,1] \times M$
then $i_1^* \beta - i_0^* \beta = \pi_* d\beta + d \pi_* \beta$
where $i_t \colon M \to M \times [0,1]$ is $i_t(m) = (m,t)$
and where $\pi_*$ is f\/iber integration.

\section{Local structure of symplectic toric manifolds}
\label{sec:background}

The purpose of this section is to set our notation
and to recall the ``local uniqueness'' result,
that symplectic toric manifolds
over the same unimodular local embedding (u.l.e.)
are locally isomorphic.  The results of the appendix
are adapted from~\cite{De}; see also~\cite{LT}.
We begin by recalling the symplectic slice representation.

\begin{Definition} \label{def:s.slice}
Let a compact Lie group~$G$ act on a symplectic manifold $(M,\omega)$
with a moment map $\mu \colon M \to {\mathfrak g}^*$.
The \emph{symplectic slice representation} at a point~$x$ of~$M$
is the linear symplectic action of the stabilizer~$G_x$
on the symplectic vector space
$V_x := T_x  (G\cdot x)^\omega /
 (T_x (G\cdot x) \cap  (T_x (G \cdot x))^\omega )$
that is induced from the linearization at $x$ of the $G_x$ action on $M$.
Here, $T_x (G\cdot x)^\omega$ denotes the symplectic perpendicular
to $T_x (G\cdot x)$ in $T_xM$.
\end{Definition}

We recall that the orbits of a Hamiltonian torus action are isotropic.
In Def\/inition~\ref{def:s.slice}, if the orbit $G \cdot x$
is isotropic, then the symplectic slice is simply
$V_x = T_x  (G\cdot x)^\omega / T_x (G\cdot x) $.

The following theorem is a consequence of the equivariant constant
rank embedding theorem of Marle~\cite{marle}; also see~\cite[Section 2]{LS}.
The proof in the case where the
group is a torus was given earlier by Guillemin and
Sternberg~\cite{GS:normal}.

\begin{Theorem}\label{thm:equiv-isotro-emb}
Let a compact Lie group $G$ act on symplectic manifolds
$(M,\omega)$ and $(M',\omega')$ with moment maps
$\mu \colon M \to {\mathfrak g}^*$ and $\mu' \colon M' \to {\mathfrak g}^*$.
Fix a point $x$ of $M$ and a point~$x'$ of~$M'$.
Suppose that $x$ and $x'$ have the same stabilizer,
their symplectic slice representations are linearly symplectically isomorphic,
and they have the same moment map value.
Then there exists an equivariant symplectomorphism
from an invariant neighbourhood of $x$ in $M$
to an invariant neighbourhood of $x'$ in $M'$
that respects the moment maps and that sends $x$ to $x'$.
\end{Theorem}

We write points in the standard torus ${\mathbb T}^\ell = {\mathbb R}^\ell/{\mathbb Z}^\ell$ as
$\ell$-tuples $[t_1, \ldots, t_\ell]$
with $(t_1, \ldots, t_\ell) \in {\mathbb R}^\ell$.
Alternatively, since ${\mathbb R}^\ell/{\mathbb Z}^\ell \simeq ({\mathbb R}/{\mathbb Z})^\ell$,
we may think of
a point on the standard $\ell$-torus ${\mathbb T}^\ell $ as an $\ell$-tuple
$(q_1,\ldots, q_\ell)$ with $q_i \in {\mathbb R}/{\mathbb Z}$.  We think of the $q_i$s as
coordinates.  Then the cotangent bundle $T^*{\mathbb T}^\ell$ has canonical
coordinates $(q_1, \ldots, q_\ell, p_1,\ldots, p_\ell)$
with $p_i\in {\mathbb R}^*$.
The symplectic form on the cotangent bundle is given by
$\omega = \sum dp_i \wedge dq_i$ in these coordinates.
The lift of the action of ${\mathbb T}^\ell$ on itself by left multiplication
to the action on
$T^*{\mathbb T}^\ell$ is Hamiltonian with an associated moment map
\begin{gather*}
T^*{\mathbb T}^\ell \to ({\mathbb R}^\ell)^* \simeq ({\mathbb R}^*)^\ell, \qquad
(q_1,\ldots, q_\ell, p_1,\ldots, p_\ell)\mapsto (p_1,\ldots,  p_\ell).
\end{gather*}

Recall the local normal form for neighbourhoods of orbits
in a symplectic toric $G$-manifold:

\begin{Lemma}\label{lem:2.6}
Let $(M,\omega, \mu \colon M \to {\mathfrak g}^*)$ be a symplectic toric $G$-manifold.
Consider a point $x$ in~$M$; denote its stabilizer by~$K$.
\begin{enumerate}[$1.$]\itemsep=0pt
\item
There exists an isomorphism $\tau_K \colon K \to {\mathbb T}^k$
such that the symplectic slice representation at~$x$
is isomorphic to the action of $K$ on ${\mathbb C}^k$ obtained from
the composition of $\tau_K$
with the standard action of ${\mathbb T}^k$ on ${\mathbb C}^k$,
which is
\begin{gather} \label{std action}
[t_1,\ldots, t_k]\cdot (z_1,\ldots, z_k) =
\big(e^{2\pi\sqrt{-1} t_1} z_1,\ldots, e^{2\pi\sqrt{-1} t_k} z_k\big).
\end{gather}

\item
Let $\tau \colon G \to {\mathbb T}^\ell \times {\mathbb T}^k$ be an isomorphism
of Lie groups such that $\tau(a) = (1,\tau_K(a))$ for all $a \in K$.
Then there exists a $G$-invariant open neighbourhood $U$
of $x$ in $M$
and a~$\tau$-equivariant open symplectic embedding
\begin{gather*}
j \colon \ U \hookrightarrow T^*{\mathbb T}^\ell \times {\mathbb C}^k
\end{gather*}
with $j(G\cdot x) = {\mathbb T}^\ell \times \{0\}$.  Here ${\mathbb T}^\ell $ acts on $
T^*{\mathbb T}^\ell$ by the lift of the left multiplication
and~${\mathbb T}^k$ acts on~${\mathbb C}^k$ by~\eqref{std action}.
Our normalization for
the symplectic form on ${\mathbb C}^k$ is
$\omega_{{\mathbb C}^k} = \frac{\sqrt{-1}}{2\pi} \sum dz_j \wedge d\overline{z}_j$,
so that
\begin{gather} \label{moment}
\mu|_U = \mu(x) + \tau^* \circ \phi \circ j,
\end{gather}
where
\begin{gather*}
\phi ((q_1, \ldots, q_\ell, p_1, \ldots p_\ell), (z_1,\ldots, z_k))
=   \Big(   (p_1,\ldots, p_\ell)   ,   \sum |z_j|^2 e_j^*   \Big)  ,
\end{gather*}
$e_1^*,\ldots,e_k^*$ is the canonical basis
of the weight lattice $({\mathbb Z}^k)^*$,
and $\tau^* \colon ({\mathbb R}^*)^\ell \times ({\mathbb R}^*)^k \to {\mathfrak g}^*$
is the isomorphism on duals of Lie algebras that is induced by $\tau$.
\end{enumerate}
\end{Lemma}

Part (1) of Lemma~\ref{lem:2.6} follows from the facts that every
linear symplectic action of a compact group preserves some compatible
Hermitian structure and that every $k$-dimensional abelian subgroup of
$U(k)$ is conjugate to the subgroup of diagonal matrices.
Part~(2) follows from Theorem~\ref{thm:equiv-isotro-emb}.

We are now ready to prove Proposition~\ref{prop:unimod-emb 1}. For the
reader's convenience we recall
its statement:

\medskip

\noindent
{\bf Proposition~\ref{prop:unimod-emb 1}.} {\em
Let $(M, \omega, \mu)$ be a symplectic toric $G$-manifold.
Then the quotient $M/G$ is naturally a manifold with corners,
and the orbital moment map $\overline{\mu} \colon M/G\to {\mathfrak g}^*$ is a u.l.e.\
$($q.v.\ Definitions~{\rm \ref{quotient}} and~{\rm \ref{def:unimodular emb})}.}

\medskip

\begin{proof}
  By Lemma~\ref{lem:2.6} we may assume 
  that $M = T^* {\mathbb T}^\ell \times {\mathbb C}^k$ with the action of $G = {\mathbb T}^\ell
  \times {\mathbb T}^k$ as in the lemma
  and that the moment map is the map $\mu \colon M \to ({\mathbb R}^*)^\ell
  \times ({\mathbb R}^*)^k$ given by
\begin{gather*}
\mu ((q_1, \ldots, q_\ell, p_1, \ldots p_\ell), (z_1,\ldots, z_k))
  =   \Big(   (p_1,\ldots, p_\ell)   ,   \sum |z_j|^2 e_j^*   \Big)   ,
\end{gather*}
where, as before, $e_1^*,\ldots,e_k^*$ is the canonical basis
of the weight lattice $({\mathbb Z}^k)^*$.
Then
\begin{gather*}
 \mu(M) = \big({\mathbb R}^\ell\big)^*   \times
           \Big\{   \sum \eta_j e_j^* \, | \,
             \eta_j \geq 0 \text{ for all $j$} \Big\}.
\end{gather*}
Hence, $\mu(M)$ is a unimodular cone; in particular, it is a manifold
with corners.
We now argue that $\mu \colon M \to \mu(M)$ is a quotient map
in the category of manifolds with corners; see Def\/inition~\ref{quotient}.
The f\/ibers of~$\mu$ are precisely the $G$-orbits.
So it remains to show that for any manifold with corners~$N$
and any $G$-invariant smooth map $f \colon M \to N$
there exists a unique smooth map
$\tilde{f} \colon \mu(M)\to N$ such that
\begin{gather*}
f = \tilde{f} \circ \mu.
\end{gather*}
Clearly, there exists a unique map $\tilde{f}$ with the above property,
and our task is to show that $\tilde{f}$ is actually smooth.
Without loss of generality we assume that $N = {\mathbb R}$.
The smoothness of~$\tilde{f}$ then follows from a special case of a~theorem of Schwarz~\cite{sc:sm}. The key point is that since the
functions $ |z_1|^2, \ldots, |z_k|^2$ generate the ring of~$ {\mathbb T}^k$
invariant polynomials on ${\mathbb C}^k$, for any smooth ${\mathbb T}^k$-invariant
function $h$ on ${\mathbb C}^k$ there is a smooth function~$\tilde{h} $ on~${\mathbb R}^k$ with $h (z_1, \ldots, z_k) = \tilde{h} (|z_1|^2,\ldots |z_k|^2)$.
\end{proof}

To f\/inish the section, we prove that symplectic toric $G$-manifolds
over the same u.l.e.\ are locally isomorphic.

\begin{Lemma} \label{lem:loc-iso}
Let $(M, \omega, \pi)$
and $(M', \omega',  \pi')$
be two symplectic toric $G$-manifolds
over the same u.l.e.\ $\psi \colon W\to {\mathfrak g}^*$.
Then for any point $w \in W$
there is a neighbourhood $U_w$ of $w$ in $W$ and an isomorphism
$\varphi \colon \pi^{-1} (U_w) \to (\pi')^{-1} (U_w)$ of symplectic
toric $G$-manifolds over $U_w \stackrel{\psi}{\to}{\mathfrak g}^* $.
\end{Lemma}

\begin{proof}
Fix $w \in W$.
Let $x$ be a point in~$\pi^{-1}(w)$.
Every invariant neighbourhood of $G \cdot x$ in~$M$
is a subset of $M$ of the form $\pi^{-1}(U_w)$
where~$U_w$ is a neighbourhood of~$w$ in~$W$.

Let $K$ be the stabilizer of $x$ in $G$, and let $k$ be the dimension
of~$K$.  Lemma~\ref{lem:2.6} implies that the symplectic slice
representation at $x$ is linearly symplectically isomorphic to the
action of~$K$ on ${\mathbb C}^k$ through an isomorphism $\tau|_K \colon K
\stackrel{\cong}{\to} (S^1)^k$.  Let $v_1^*,\ldots,v_k^*$ denote the
basis of the weight lattice ${\mathbb Z}_K^*$ that corresponds under~$\tau|_K$
to the standard basis of the weight lattice~$({\mathbb Z}^k)^*$.  Thus,
$v_1^*,\ldots,v_k^*$ represent the weights of the~$K$ action on~${\mathbb C}^k$.  Let $v_1,\ldots,v_k$ be the dual basis, in~${\mathbb Z}_K$, and let
$\epsilon = \mu(x)$.  The equation~\eqref{moment} for the moment map
implies that every neighbourhood of~$x$ contains a smaller invariant
neighbourhood of~$x$ whose moment map image is a neighbourhood of~$\epsilon$ in the cone $C_{\{v_1,\ldots,v_k\},\epsilon}$ (cf.\
Def\/inition~\ref{def:unimod_cone}).

By combining the above discussion with  Lemma~\ref{lem:prim-norm} we see that
the symplectic slice representation is determined up to linear
symplectic isomorphism by the image of
an arbitrary suf\/f\/iciently small invariant neighbourhood of $G \cdot x$.
This image is exactly $\psi(U_w)$, where~$U_w$ is an arbitrary
suf\/f\/iciently small
neighbourhood of~$w$ in~$W$.  So the germ of $\psi$ at $w$ determines
the symplectic slice representation up to linear symplectic
isomorphism.  Clearly, this germ also determines the moment map value,~$\mu(x)$, which is equal to~$\psi(w)$.  The result then follows from
Theorem~\ref{thm:equiv-isotro-emb}.
\end{proof}

\subsection*{Acknowledgments}
This research is partially supported by the Natural Sciences
and Engineering Research Council of Canada  and by the National
Science Foundation.

We are grateful to Chris Woodward and to Daniele Sepe
for valuable discussions.
We are also grateful to our anonymous referees for their valuable comments.

\pdfbookmark[1]{References}{ref}
\LastPageEnding

\end{document}